# SUBSPACES OF SEPARABLE $L_1$-PREDUALS: $W_\alpha$ EVERYWHERE

EMANUELE CASINI, ENRICO MIGLIERINA, AND ŁUKASZ PIASECKI

ABSTRACT. The spaces $W_\alpha$ are the Banach spaces whose duals are isometric to $\ell_1$ and such that the standard basis of $\ell_1$ is $w^*$-convergent to $\alpha \in \ell_1$. The core result of our paper proves that an $\ell_1$-predual $X$ contains isometric copies of all $W_\alpha$, where the norm of $\alpha$ is controlled by the supremum of the norms of the $w^*$-cluster points of the extreme points of the closed unit ball in $\ell_1$. More precisely, for every $\ell_1$-predual $X$ we have

$$r^*(X) = \sup\left\{\|g^*\| : g^* \in (\text{ext } B_{\ell_1})'\right\} = \sup\left\{\|\alpha\| : \alpha \in B_{\ell_1},\ W_\alpha \subset X\right\}.$$

We also prove that, for any $\varepsilon > 0$, $X$ contains an isometric copy of some space $W_\alpha$ with $\|\alpha\| > r^*(X) - \varepsilon$ which is $(1 + \varepsilon)$-complemented in $X$. From these results we obtain several consequences. First we provide a new characterization of separable $L_1$-preduals containing an isometric copy of a space of affine continuous functions on a Choquet simplex. Then, we prove that an $\ell_1$-predual $X$ contains almost isometric copies of the space $c$ of convergent sequences if and only if $X^*$ lacks the stable $w^*$-fixed point property for nonexpansive mappings.

## 1. INTRODUCTION

A Banach space $X$ such that its dual $X^*$ is isometric to $L_1(\mu)$ for some measure $\mu$ is usually called an $L_1$-predual (or a Lindenstrauss space). This class of spaces received many attentions in the '60 and in the '70 of the last century and many results had been obtained (see the monograph [23] and the reference therein). One of the most important result about this topic has been obtained by Zippin [37]: he proved that an isometric copy of the Banach space $c_0$ of sequences converging to 0 is contained in every separable $L_1$-predual. He also showed that there always exists a copy of $c_0$ which is 1-complemented. Moreover, we know from Lazar and Lindenstrauss [25, 26] that the space $C(\Delta)$ of continuous functions on the Cantor set $\Delta$ is contained in every separable $L_1$-predual with nonseparable dual. Therefore, all such spaces contain an isometric copy of the Banach space $c$ of convergent sequences.

The spaces $c_0$ and $c$ share two common features: first, both of them are preduals of $\ell_1$ and, second, they induce a $w^*$-topology on $\ell_1$ such that the standard basis of $\ell_1$ is $w^*$-convergent. Given $\alpha$ in the closed unit ball $B_{\ell_1}$ of $\ell_1$, we denote by $W_\alpha$ the $\ell_1$-predual such that the standard basis of $\ell_1$ is $w^*$-convergent to $\alpha$. In [2], we provided a concrete representation of such spaces as hyperplanes of $c$, namely:

$$W_\alpha = \left\{ x \in c : \lim_{i \to \infty} x(i) = \sum_{i=1}^{\infty} \alpha(i) x(i) \right\}.$$

To our best knowledge, the spaces $W_\alpha$ appeared only as examples in [27], p. 78, in [24], p. 363 and in [19], p. 239. The first systematic study of their properties was made by the present authors in [2, 3]. These spaces found applications in the study of $\sigma(X^*, X)$-fixed point property in the dual of an $L_1$-predual $X$, where $\sigma(X^*, X)$ denotes the weak* topology on $X^*$ induced by the $L_1$-predual $X$ (see [6, 7, 3, 4, 34, 33, 13, 35] and see also [17, 18, 9]). Moreover, in [8] it turned out that one of these spaces provides a surprising counterexample to a result about the presence of an isometric copy of $c$ in every $L_1$-predual whose closed unit ball has an extreme point (see Section 4 in [37]). Recently, it occured that some of $W_\alpha$ spaces play a crucial role in the description of all $L_1$-preduals failing the extension of compact operators property (see [8] and Theorem 3.4 in [13]). Topological and metric properties of the space $(\mathcal{H}, d)$, where $\mathcal{H}$ is the set of all $W_\alpha$ spaces, $d$ denotes the Banach–Mazur distance and we identify hyperplanes that are Banach–Mazur distance 0 apart, are studied in details in (see [14, 33]). Finally, spaces $W_\alpha$ give interesting examples in some other topics (see [31]).

In the recent paper [38], taking into account [8], Zippin proved that each separable $L_1$-predual contains $c$ or an $\ell_1$-predual hyperplane of $c$. One of the main aims of the present paper is to provide a deeper description of the subspaces contained in an $\ell_1$-predual $X$, by providing a "quantitative" version of Zippin's result. By means of the class of $W_\alpha$ spaces, this approach makes clear which hyperplanes of $c$ are contained in $X$. Indeed, we show how $W_\alpha$ spaces can be seen as "constitutive bricks" of $\ell_1$-preduals because they are ubiquitous and simple, but remain $\ell_1$-preduals in themselves.

The paper is organized as follows.

In Section 2 we collect some preliminaries. In particular, we briefly describe the structure of separable $L_1$-preduals made in terms of finite-dimensional $\ell_\infty^n$ spaces ([25]), that we call "cubic" structure of an $\ell_1$-predual. Here, $\ell_\infty^n$ denotes the space $\mathbb{R}^n$ endowed with the sup norm.

In Section 3 we present an extension of an old result about extreme points of the closed unit ball of the dual space of a subspace of the space $C(K)$ of continuous functions on a Hausdorff compact $K$ (see Lemma V.8.6 in [10]) to the dual of an $\ell_1$-predual (see Lemma 3.1 and Remark 3.2). It will play an important role in several points of our study.

Section 4 establishes a new characterization of the existence of a 1-complemented isometric copy of $c$ in a separable $L_1$-predual $X$. This condition involves the behaviour of the extreme points of a face of the dual ball with respect to the $\sigma(X^*, X)$-topology. One of the sides of the proof is based on a powerful result proved by Gasparis [12].

In Section 5 we further develop the theory of $W_\alpha$ spaces in order to get some results that will be useful in the sequel. After a brief summary of known properties, we study their representations by matrix that sum up their cubic structure. Moreover, we carefully study the problem to isometrically embed a space $W_\beta$ in a space $W_\alpha$, where $\alpha, \beta \in B_{\ell_1}$.

Sections 6 and 7 are the core of our paper. In all these sections, our starting point are assumptions about the set of $w^*$-cluster points of the extreme points of the closed unit ball of $X^*$. To be more precise, we consider the radius of the smallest ball containing the set $(\text{ext} B_{X^*})'$ of the $\sigma(X^*, X)$-cluster points of the extreme points of $B_{X^*}$, i.e.,

$$r^*(X) = \sup\{\|g^*\| : g^* \in (\text{ext } B_{X^*})'\},$$

without assuming any information about $X$ itself. Indeed, in Section 6, we prove that, if $X$ is an $\ell_1$-predual, it holds

$$r^*(X) = \sup\{\|\alpha\| : \alpha \in B_{\ell_1}, W_\alpha \subset X\}.$$

This result has some interesting consequences; for example, it allows us to obtain a new characterization of separable $L_1$-preduals $X$ containing an isometric copy of a space of affine continuous functions on a Choquet simplex. Section 7 contains a result showing that every $\ell_1$-predual $X$ contains a $(1 + \varepsilon)$-complemented isometric copy of $W_\alpha$ with $\alpha = \left(\frac{r^*(X)}{1+\varepsilon}, 0, 0, \dots\right)$, for every $\varepsilon > 0$. In addition, by combining this fact with the results of Subsection 5.2, we show that, for every $\alpha \in \ell_1$ with $\|\alpha\| < r^*(X)$, $X$ contains a complemented subspace which is isometric to $W_\alpha$.

In Section 8, we show that the condition $r^*(X) = 1$ is equivalent to the existence of almost isometric copies of $c$ in $X$. It is worth pointing out the main reason of interest of this section: namely, by Theorem 3.5 in [7], we prove that the existence of almost isometric copies of $c$ in $X$ is equivalent to the lackness of the stable $\sigma(X^*, X)$-fixed point property for nonexpansive mappings. Thus, we obtain a link between two, relatively far, parts of Banach spaces theory.

We conclude the paper with Section 9, where we show that for almost all classes of $L_1$-preduals introduced in literature (see [29, 22, 23] and the references therein), it is possible to find a concrete example of a space belonging to each class by using $W_\alpha$ spaces or concrete representations of $\ell_1$-preduals introduced in [5].

## 2. Preliminaries

Let $B_X$ and $S_X$ denote, respectively, the closed unit ball and the unit sphere in a real Banach space $X$, and $X^*$ denote the topological dual of $X$. If $K$ is a nonempty, convex subset of a linear topological space, then by ext $K$ we denote the set of all extreme points of $K$. A convex subset $F$ of $B_X$ is called a *face* of $B_X$ if for every $x, y \in B_X$ and $\lambda \in (0, 1)$ such that $(1 - \lambda)x + \lambda y \in F$ we have $x, y \in F$. A face $F$ of $B_X$ is named a *proper face* if $F \neq B_X$. Let $A$ be a subset of $X$, we denote by $\overline{A}$, span($A$) and by conv($A$) the norm closure of $A$, the linear span of $A$ and the convex hull of $A$, respectively. Moreover, if $A \subseteq X^*$, we denote by $\overline{A}^{\sigma(X^*, X)}$ the $\sigma(X^*, X)$-closure of $A$ and by $A'$ the set of $\sigma(X^*, X)$-cluster points of $A$, where $\sigma(X^*, X)$ denotes the usual $w^*$-topology on $X^*$ induced by $X$. Finally, the standard basis of $\ell_1$ is denoted by $\{e_n^*\}_{n=1}^\infty$.

It is well known that $c^*$ can be isometrically identified with $\ell_1$ in this way: for every $x^* \in c^*$ there exists a unique $f^* = (f^*(1), f^*(2), \dots) \in \ell_1$ such that

$$x^*(x) = f^*(1) \lim_{i \to \infty} x(i) + \sum_{i=1}^\infty f^*(i+1)x(i) = f^*(x),$$

where $x = (x(1), x(2), \dots) \in c$. It is easily seen that the the elements $e_n^*$ of the standard basis of $\ell_1$ are functionals on $c$ such that $e_n^*(x) = x(n-1)$ for every $n \geq 2$ and for every $x \in c$. Moreover, it holds

$$e_1^*(x) = \lim_{n \to \infty} x(n) \quad \text{and} \quad e_n^* \xrightarrow[n \to \infty]{\sigma(\ell_1, c)} e_1^*.$$

We recall that an $L_1$-predual space $X$ is isometric to the space $A(S)$ of affine continuous function on a Choquet simplex $S$ if and only if $B_X$ has at least one extreme point (see [36]). Moreover, a Banach space $X$ is polyhedral if the unit balls of all its finite-dimensional subspaces are polytopes (see [21]).

Finally, in $\ell_1$ we introduce the following notation: let $x^* \in \ell_1$, we put

$$\text{supp}(x^*) := \{i \in \mathbb{N} : x^*(i) \neq 0\}.$$

2.1. **"Cubic" structure of separable $L_1$-preduals.** Here we recall a, now classical, result about the structure of separable $L_1$-preduals. This result was originally proved in [25]; here we follow the presentation given by Zippin in the introduction of the paper [38].

Let $X$ be a separable $L_1$-predual, then there exists a normalized monotone basis $\{x_i\}_{i=1}^\infty$ of $X$ such that $X_n = \text{span}(\{x_i\}_{i=1}^n)$ is isometric to $\ell_\infty^n$ and the sequence $\{e_i^n\}_{i=1}^n$ is a basis for $X_n$ satisfying the following properties:

$$(2.1) \qquad\qquad e_n^n = x_n \quad \text{for every } n \in \mathbb{N},$$

$$(2.2) \qquad\qquad \left\| \sum_{i=1}^n c_i e_i^n \right\| = \max_{1 \le i \le n} |c_i| \quad \text{for all } \{c_i\}_{i=1}^n \subset \mathbb{R}$$

and

$$(2.3) \qquad \text{for every } 1 \le i \le n, \quad e_i^n = e_i^{n+1} + a_i^n e_{n+1}^{n+1} \quad \text{where } \sum_{i=1}^n |a_i^n| \le 1.$$

It is worth pointing out that also the reverse implication holds (see Theorem 2 in [25]). By using these construction, it is clear that a separable $L_1$-predual can be represented by an infinite triangular matrix $A = (a_i^n)$ (see [26, 30]). Now we introduce the functionals $\{\phi_j\}_{j=1}^\infty$ defined on $\bigcup_{n=1}^\infty X_n$ by

$$(2.4) \qquad \phi_j\left( \sum_{i=1}^n c_i e_i^n \right) = c_j \quad \text{for all } j \le n \text{ and } \{c_i\}_{i=1}^n \subset \mathbb{R}.$$

All these functionals are well-defined and each of them can be uniquely extended to a a linear functional on $X = \overline{\bigcup_{n=1}^\infty X_n}$. We also point out that the sequence $\{\phi_j\}_{j=1}^\infty$ is isometrically equivalent to the standard basis of $\ell_1$ (i.e., $\left\| \sum_{j=1}^n a_j \phi_j \right\| = \sum_{j=1}^n |a_j|$ for all $n \in \mathbb{N}$ and sequences $\{a_j\}_{j=1}^\infty$). Moreover, the following result holds.

*Lemma* 2.1. Let $\{\phi_j\}_{j=1}^\infty$ be the set of functionals defined as above. Then, the two following properties hold:

(i) $\{\phi_j\}_{j=1}^\infty \subseteq \text{ext}\, B_{X^*}$;

(ii) $\overline{\{\pm\phi_j\}_{j=1}^\infty}^{\sigma(X^*,X)} = \overline{\text{ext}\, B_{X^*}}^{\sigma(X^*,X)}$.

The proof of (i) is the proof of Lemma 1 in [37], whereas the proof of (ii) is the proof of Lemma 1.2 in [30].

The following remark comes from [32] and it plays a key role in one of our main results.

*Remark* 2.2. The first element $x_1$ of the monotone basis $\{x_i\}_{i=1}^{\infty}$ of $X$ can be chosen to be an arbitrary vector of norm one in $X$.

The principal interest of the following result is that it allows one to build up a subspace $U$ of a separable $L_1$-predual such that $U$ itself is an $L_1$-predual.

*Lemma* 2.3. ([38], Lemma 1) Let $\{n(k)\}_{k=1}^{\infty}$ be a strictly increasing sequence of the integers such that $n(1) = 1$. Then the sequence $\left\{e_{n(k)}^{n(k)}\right\}_{k=1}^{\infty}$ is a monotone basis for the subspace $U = \overline{\mathrm{span}\left(\left\{e_{n(k)}^{n(k)}\right\}_{k=1}^{\infty}\right)} \subseteq X$. Moreover, in each subspace

$$U_m = \mathrm{span}\left(\left\{e_{n(k)}^{n(k)}\right\}_{k=1}^{m}\right)$$

there is a basis $\left\{u_i^m\right\}_{i=1}^{m}$ satisfying the following three conditions

$$(2.5) \qquad u_m^m = e_{n(m)}^{n(m)} \quad \text{for every } m \in \mathbb{N},$$

$$(2.6) \qquad \left\|\sum_{i=1}^{m} c_i u_i^m\right\| = \max_{1 \le i \le m} |c_i| \quad \text{for all } \{c_i\}_{i=1}^{m} \subset \mathbb{R}$$

and

$$(2.7) \qquad \text{for every } 1 \le i \le m, \quad u_i^m = u_i^{m+1} + b_i^m e_{m+1}^{m+1}, \quad \text{where } \sum_{i=1}^{m} \left|b_i^m\right| \le 1.$$

Finally, let $\psi_j = \phi_{n(j)}|_U$; then $\psi_j\left(u_i^m\right) = \delta_{i,j}$ for all $1 \le i, j \le m$ and, for every $u \in U$, we have

$$(2.8) \qquad \|u\| = \sup_{j \ge 1} \left|\psi_j(u)\right| = \sup_{j \ge 1} \left|\phi_{n(j)}(u)\right|.$$

## 3. Extreme points of the dual of a subspace of an $\ell_1$-predual

In this brief section we prove a result that describes all the extreme points of the dual ball of a subspace $X$ of a given $\ell_1$-predual $Y$. This result is an extension, in the separable case, of Lemma V.8.6 in [10], where the space $Y$ is the space $C(Q)$ of continuous functions defined on the compact Hausdorff set $Q$.

*Lemma* 3.1. Let $X$ be a closed subspace of an $\ell_1$-predual $Y$ and let us denote by $\{e_n^*\}_{n=1}^{\infty}$ the standard basis of $Y^*$. Then, there exists a function $k : \mathrm{ext}\, B_{X^*} \to \mathbb{N}$ such that for every $f^* \in \mathrm{ext}\, B_{X^*}$ it holds

$$(3.1) \qquad f^* = \delta(f^*)\, e_{k(f^*)}^*\big|_X,$$

where $\delta(f^*)$ is a real number such that $|\delta(f^*)| = 1$ for every $f^* \in \mathrm{ext}\, B_{X^*}$.

*Proof.* Let us consider the set

$$A = \left\{ \pm e_n^* \big|_X \right\}_{n=1}^{\infty}.$$

By contradiction, let us suppose that there exists $g^* \in \operatorname{ext} B_{X^*} \setminus A$. Then there exists $h^* \in Y^*$ such that $h^*|_X = g^*$ and $\|h^*\| = 1$. Since $g^* \notin A$, it holds

$$(3.2) \qquad h^* = \sum_{j \in \Delta} h^*(j) e_j^*,$$

where $\Delta = \{ k \in \mathbb{N} : h^*(k) \neq 0 \}$. Since

$$1 = \|g^*\| = \left\| \sum_{j \in \Delta} h^*(j) \, e_j^* \big|_X \right\| \leq \sum_{j \in \Delta} |h^*(j)| \, \big\| e_j^* \big|_X \big\| \leq \sum_{j \in \Delta} |h^*(j)| \, \big\| e_j^* \big\| = 1,$$

we obtain that $\big\| e_j^*|_X \big\| = 1$ for $j \in \Delta$. It is easy to see that $\Delta$ has at least two elements. Indeed, if $\Delta$ is a singleton then $g^* \in A$, a contradiction. For $k \in \Delta$ let $u^* = h^*(k) e_k^*|_X$ and $v^* = \sum_{j \in \Delta, \, j \neq k} h^*(j) e_j^*|_X$. Then, $g^* = u^* + v^*$ and $u^* \neq 0$ imply $v^* \neq 0$ since if $v^* = 0$ then, $g = h^*(k) e_k^*|_X = \pm e_k^*|_X$, that is, $g^* \in A$. Moreover, if there exists $\lambda$ such that $v^* = \lambda u^*$, we obtain $g^* = (1 + \lambda) h^*(k) e_k^*|_X = \pm e_k^*|_X$, contrary to $g^* \in A$. Therefore, we conclude that $v^* \neq \lambda u^*$ for every $\lambda$. Finally, it holds:

$$1 = \|g^*\| = \|u^* + v^*\| \leq \|u^*\| + \|v^*\| \leq |h^*(k)| + \sum_{j \in \Delta, \, j \neq k} |h^*(j)| = 1.$$

Hence, we get:

$$g^* = u^* + v^* = \|u^*\| \frac{u^*}{\|u^*\|} + (1 - \|u^*\|) \frac{v^*}{\|v^*\|}.$$

Since $g^* \in \operatorname{ext} B_{X^*}$ it follows that $\frac{u^*}{\|u^*\|} = \frac{v^*}{\|v^*\|}$, which is impossible because $v^*$ is not a multiple of $u^*$. $\qquad \square$

*Remark* 3.2. As pointed out in a private communication by Carlo A. De Bernardi to the authors, the previous result holds in a very general framework. Indeed, if $X$ is a closed subspace of a dual $Y^*$ of a normed space $Y$, then by combining a consequence of the well-known Krein–Milman theorem (see Problem C, p. 132 in [20]) and properties of weak* topology, we can prove that each extreme point of the closed unit ball of $X$ can be seen as a restriction of an extreme point of $B_{Y^*}$. Nevertheless, here we prefer to present a direct version of the proof of Lemma 3.1, since it is very basic and it does not use Krein–Milman theorem.

## 4. Existence of a subspace isometric to the space $c$

The main result of this section characterizes the $\ell_1$-preduals containing an isometric copy of the space $c$. The presence of an isometric copy of $c$ in an $\ell_1$-predual $X$ is equivalent both to polyhedrality of $X$ and to the property: there exists $x \in S_X$ such that

$$(4.1) \qquad \sup \{ x^*(x) : x^* \in \operatorname{ext} (B_{X^*}) \setminus D(x) \} = 1,$$

where $D(x) = \{ x^* \in S_{X^*} : x^*(x) = 1 \}$. This result was proved in Theorem 4.1 of [8]. In the same paper, a necessary condition for the presence of an isometric copy of $c$ in a generic separable Banach space was found (Theorem 2.1 in [8]).

By restricting this condition to $\ell_1$-preduals, the following necessary condition can be obtained:

*Corollary* 4.1. (Corollary 2.2 in [8]) Let $X$ be a predual of $\ell_1$. If $X$ contains a subspace isometric to $c$, then there exist $x \in B_X$ and a subsequence $\{e^*_{n_k}\}_{k \in \mathbb{N}}$ of the standard basis $\{e^*_n\}_{n \in \mathbb{N}}$ such that

(1) $e^*_{n_k} \xrightarrow[k \to \infty]{\sigma(\ell_1, X)} e^*$,

(2) $\text{supp}(e^*_{n_k}) \cap \text{supp}(e^*) = \emptyset$ for every $k \in \mathbb{N}$,

(3) $e^*_{n_k}(x) = e^*(x) = 1$ for every $k \in \mathbb{N}$.

The following example shows that Corollary 4.1 gives only a necessary condition for the existence of a subspace isometric to $c$.

*Example* 4.2. In this example we represent the ordinal interval $[1, \omega 2]$ (with the order topology) by the subset of $\mathbb{R}^2$ (with the Euclidean topology) defined as

$$\Omega^2 = \Omega_0^2 \cup \left\{ \left( \frac{1}{i}, 0 \right) \in \mathbb{R}^2 : i = 1, 2 \right\},$$

where

$$\Omega_0^2 = \left\{ \left( \frac{1}{i}, \frac{1}{j} \right) \in \mathbb{R}^2 : i = 1, 2 \text{ and } j = 1, 2, \dots \right\}.$$

Let $x_1^*, x_2^* \in \ell_1(\Omega_0^2)$ be defined by

$$x_1^*(i, j) = \begin{cases} 0 & \text{for } i = 1 \text{ and } j \geq 1 \\ \frac{1}{2^j} & \text{for } i = 2 \text{ and } j \geq 1 \end{cases},$$

and

$$x_2^*(i, j) = \begin{cases} \frac{1}{2^j} & \text{for } i = 1 \text{ and } j \geq 1 \\ 0 & \text{for } i = 2 \text{ and } j \geq 1 \end{cases}.$$

Then, we consider the space

$$W_{x_1^*, x_2^*} = \left\{ f \in C(\Omega^2) : f(1, 0) = \sum_{j=1}^\infty \frac{1}{2^j} f\left( \frac{1}{2}, \frac{1}{j} \right) \text{ and } f\left( \frac{1}{2}, 0 \right) = \sum_{j=1}^\infty \frac{1}{2^j} f\left( 1, \frac{1}{j} \right) \right\}.$$

By Theorem 2.1 in [5], $W^*_{x_1^*, x_2^*} = \ell_1(\Omega_0^2) = \ell_1$ and

$$e^*_{1,j} \xrightarrow{\sigma(\ell_1(\Omega_0^2), W_{x_1^*, x_2^*})} x_1^* = \left( 0, 0, 0, \dots, \frac{1}{2}, \frac{1}{4}, \frac{1}{8}, \dots \right),$$

and

$$e^*_{2,j} \xrightarrow{\sigma(\ell_1(\Omega_0^2), W_{x_1^*, x_2^*})} x_2^* = \left( \frac{1}{2}, \frac{1}{4}, \frac{1}{8}, \dots, 0, 0, 0, \dots \right),$$

where $\left\{ e^*_{1,j}, e^*_{2,j} \right\}_{j=1}^\infty$ is the standard basis of $\ell_1(\Omega_0^2)$.

Let us consider the element $\bar{x} \in S_{W_{x_1^*, x_2^*}}$ defined by

$$\bar{x}\left( 1, \frac{1}{j} \right) = \bar{x}(1, 0) = \bar{x}\left( \frac{1}{2}, \frac{1}{j} \right) = \bar{x}\left( \frac{1}{2}, 0 \right) = 1 \text{ for every } j \geq 1.$$

It is easy to see that $e^*_{1,j}(\bar{x}) = x_1^*(\bar{x}) = 1$ and $\text{supp}(e^*_{1,j}) \cap \text{supp}(x_1^*) = \emptyset$ for every $j \geq 1$. Therefore, conditions (1), (2), and (3) of Corollary 4.1 are satisfied.

Nevertheless, the space $W_{x_1^*, x_2^*}$ does not contain an isometric copy of $c$. Indeed, the only element $x \in S_{W_{x_1^*, x_2^*}}$ such that

$$x(1, 0) = x\left(\frac{1}{2}, 0\right) = 1$$

is the element $\bar{x}$ considered above. Hence, it is easy to see that there is no element in $S_{W_{x_1^*, x_2^*}}$ satysfying condition (4.1).

Our aim is to strengthen condition (2) in Corollary 4.1, in order to obtain a necessary and sufficient condition for the existence of an isometric copy of $c$ in an $\ell_1$-predual.

In the proof of our result we need some properties of the set $D(x)$. Namely, we recall that $D(x)$ is a $w^*$-compact proper face for every $x \in S_X$ and consequently $\operatorname{ext} D(x) = D(x) \cap \operatorname{ext} B_{X^*} \neq \emptyset$ by the Krein–Milman theorem. Moreover, for every $x \in S_{\ell_1}$, the set $D(x)$ can be described as follows:

$$D(x) = \left\{ \sum_{i \in \Lambda} \delta_i d_i^* : \sum_{i \in \Lambda} \delta_i = 1, \ \delta_i \geq 0 \text{ for all } i \in \Lambda \right\},$$

where $\left\{ d_i^* \right\}_{i \in \Lambda} = \operatorname{ext} D(x)$ and $\Lambda \subset \mathbb{N}$ is a set of indexes.

*Theorem* 4.3. Let $X$ be a predual of $\ell_1$. Then, $X$ contains an isometric copy of $c$ if and only there exist $\bar{x} \in S_X$ and a sequence $\{v_n^*\}_{n=1}^\infty \subseteq \operatorname{ext} B_{X^*}$ such that

(i) $v_n^* \xrightarrow[n \to \infty]{\sigma(X^*, X)} \bar{x}^*$;

(ii) $v_n^*(\bar{x}) = \bar{x}^*(\bar{x}) = 1$ for every $n \in \mathbb{N}$;

(iii) $\operatorname{supp}(v_n^*) \cap \operatorname{supp}(d^*) = \emptyset$ for every $n \in \mathbb{N}$ and for all $d^* \in (\operatorname{ext} D(\bar{x}))'$.

*Proof.* Let $Z$ be a closed subspace of $X$ isometric to $c$. In the sequel of the proof we identify $Z$ with $c$ for the sake of simplicity. Let $\{e_n^*\}_{n=1}^\infty$ be the standard basis of $c^* = \ell_1$. For each $n$, we consider a norm preserving extension of $e_n^*$ to the whole $X$; we denote it by $x_n^*$. By Lemma 3.1, we can always choose $x_n^*$ such that $\{x_n^*\}_{n=1}^\infty \subseteq \operatorname{ext} B_{X^*}$. Now, we consider a subsequence $\left\{ x_{n_j}^* \right\}_{j=1}^\infty$ of $\{x_n^*\}_{n=1}^\infty$ such that $n_1 > 1$ and

$$x_{n_j}^* \xrightarrow[j \to \infty]{\sigma(X^*, X)} \bar{x}^*,$$

where $\bar{x}^* \in B_{X^*}$. Since $x_{n_j}^* \in B_{X^*} \cap D(\bar{x})$ for all $j$, we have

$$(4.2) \qquad \operatorname{ext} D(\bar{x}) = \left\{ x_{n_j}^* \right\}_{j=1}^\infty \cup \{h_k^*\}_{k \in \Lambda_1},$$

where $\Lambda_1 \subseteq \Lambda$ and $\operatorname{supp}(x_{n_j}^*) \cap \operatorname{supp}(h_k^*) = \emptyset$ for all $j = 1, 2, \ldots$ and for all $k \in \Lambda_1$. As we did in the proof of Proposition 3.1 in [3], we consider the subspace of $c$ (that is isometric to $c$ in its turn)

$$Y = \left\{ y \in c : \lim y(n) = y(s) \text{ for all } s \in \mathbb{N} \setminus \{n_j - 1\}_{j=1}^\infty \right\}$$

and the norm one projection $P : X \to Y$ defined by:

$$P(x) = \bar{x}^*(x)\bar{x} + \sum_{j=1}^\infty \left( x_{n_j}^* - \bar{x}^* \right)(x) e_{n_j - 1},$$

where $e_{n_j-1}$ is the element of $c$ with the $\left(n_j - 1\right)$-th component equal to 1 whereas all the others are 0. Now, let $d^* \in (\text{ext } D(\bar{x}))'$. Then, it holds

$$d^* = \sum_{j=1}^{\infty} \delta_{n_j} x_{n_j}^* + \sum_{k \in \Lambda_1} \delta_k h_k^*,$$

where $\delta_{n_j} \geq 0$ for all $j$, $\delta_k \geq 0$ for all $k \in \Lambda_1$, and

$$(4.3) \qquad\qquad \sum_{j=1}^{\infty} \delta_{n_j} + \sum_{k \in \Lambda_1} \delta_k = 1.$$

Since $d^* \in D(\bar{x})$ and $\bar{x} \in Y$, we have $d^*(\bar{x}) = 1 = d^*|_Y(\bar{x})$. Now, we consider the adjoint $P^*$ of the projection $P$. Hence, for every $y^* \in Y^*$ and $x \in X$, we get

$$P^*(y^*)(x) = y^*(1)\bar{x}^*(x) + \sum_{j=1}^{\infty} y^*(n_j + 1)\left(x_{n_j}^* - \bar{x}^*\right)(x) + \sum_{k \in \mathbb{N} \setminus \{n_j\}_{j=1}^{\infty}} y^*(k)\bar{x}^*(x).$$

Then, by recalling that $P(\bar{x}) = \bar{x}$, we obtain

$$1 = d^*|_Y(\bar{x}) = d^*\left(P(\bar{x})\right) = P^*\left(d^*|_Y\right)(\bar{x}) \leq \sum_{k \in \Lambda_1} \delta_k \leq 1.$$

By (4.3), we conclude that $\sum_{j=1}^{\infty} \delta_{n_j} = 0$. Hence, $\text{supp}(d^*) \cap \text{supp}(x_{n_j}^*) = \emptyset$ for all $j$. The assertions (i), (ii) and (iii) hold for $v_j^* = x_{n_j}^*$ for every $j$.

Now we pass to prove the reverse implication. First, we note that (ii) implies that $\{v_n^*\}_{n=1}^{\infty}$ and $v^*$ are contained in $D(x)$. Let us consider the set $K = \text{ext } D(\bar{x})$. We see at once that:

- $K$ is a countable subset of $\text{ext } B_{X^*}$;
- $K \cap (-K) = \emptyset$ because, if $K$ contains both the extreme points $h^*$ and $-h^*$, then $0 \in D(\bar{x})$, which is impossible;
- the set $U = \overline{\text{span}(K)}$ is $\sigma(\ell_1, X)$-closed in $X^*$. Indeed, the set $K$ is a sequence isometrically equivalent to the standard $\ell_1$-basis. By applying the Choquet representation and the Krein–Milman theorem we obtain that $B_U = \text{conv}\left(D(\bar{x}) \cup (-D(\bar{x}))\right) = \overline{\text{conv}\left(K \cup (-K)\right)}^{\sigma(\ell_1, X)}$. The assertion now follows from the Krein–Smulian theorem (Theorem V.5.7 in [10]). This fact follows also from Lemma 1 in [1].

We are now in position to apply Theorem 1.1 in [12]. By using this result we obtain that there exists a $\sigma(\ell_1, X)$-continuous norm one projection $P_1$ from $X^*$ onto $U$. Now the task is to find a $\sigma(\ell_1, X)$-continuous norm one projection $P_2$ from $U$ onto a subspace $Q$ of $U$ such that $Q$ is isometric to $c^*$ (endowed with its natural $\sigma(\ell_1, c)$-topology). In order to define such a projection, first we split the countable set $\mathcal{I}$ of indexes that enumerate the elements of $K$ into two subsets. Let $\mathcal{I}_1$ be the set of indexes enumerating the elements of $\{v_n^*\}_{n=1}^{\infty}$ and $\mathcal{I}_2 = \mathcal{I} \setminus \mathcal{I}_1$. If we denote by $g_k^*$ the elements such that the index $k \in \mathcal{I}_2$, then we can describe the subspace $U$ as follows:

$$U = \left\{ u^* \in X^* : u^* = \sum_{n=1}^{\infty} \eta_n v_n^* + \sum_{k \in \mathcal{I}_2} \gamma_k g_k^*, \ \sum_{n \in} |\eta_n| < \infty, \ \sum_{k \in \mathcal{I}_\epsilon} |\gamma_k| < \infty \right\}.$$

Let us define the linear map $P_2 : U \to X^*$ by

$$P_2(u^*) = \sum_{n=1}^{\infty} \eta_n v_n^* + \left( \sum_{k \in \mathcal{I}_2} \gamma_k \right) \bar{x}^*.$$

We remark that:

- $P_2$ is a projection from $U$ onto its closed subspace $Q$ defined by

$$Q = \left\{ q^* \in X^* : q^* = \sum_{n=1}^{\infty} \eta_n v_n^* + \theta \bar{x}^* \sum_{k=1}^{\infty} |\eta_k| < \infty, \ \theta \in \mathbb{R} \right\}.$$

  Indeed, $\bar{x}^* \in D(\bar{x})$ and by the assumptions (i) and (iii) we obtain $\text{supp}(\bar{x}^*) \subseteq \mathcal{I}_2$. Hence, there exists $\{\varepsilon_i\}_{i \in \mathcal{I}_2} \subset [0,1]$ such that $\sum_{i \in \mathcal{I}_2} \varepsilon_i = 1$ and $\sum_{i \in \mathcal{I}_2} \varepsilon_i g_i^* = \bar{x}^*$. Let us consider $q = \sum_{n=1}^{+\infty} \eta_n v_n^* + \theta \bar{x}^* \in V$, then

$$P_2(q) = P_2 \left( \sum_{n=1}^{+\infty} \eta_n v_n^* + \theta \left( \sum_{i \in \mathcal{I}_\in} \varepsilon_i g_i^* \right) \right) = \sum_{n=1}^{+\infty} \eta_n v_n^* + \theta \left( \sum_{i \in \mathcal{I}_\in} \varepsilon_i \right) v^* = q.$$

- $V$ is a $\sigma(\ell_1, X)$-closed subspace of $X^*$ (to prove this fact it is sufficient to remember that $\sigma(\ell_1, X) - \lim_{n \to \infty} v_n^* = \bar{x}^*$ by assumption (i));
- $P_2$ is a $\sigma(\ell_1, X)$-continuous map by assumptions (i), (iii) and Lemma 2 in [1];
- $\|P_2\| = 1$ because $\mathcal{I}_1 \cap \mathcal{I}_2 = \emptyset$.

By combining the two projections $P_1$ and $P_2$ introduced above we now define the linear map $P : X^* \longrightarrow Q$ such that $P(x^*) = P_2(P_1(x^*))$ for every $x^* \in X^*$. The map $P$ is a $\sigma(X^*, X)$-continuous and norm one projection from $X^*$ onto $Q$.

Let $\{f_n^*\}_{n=1}^{+\infty}$ denote the standard basis of $c^* = \ell_1$. Then the map, defined by

$$T(\bar{x}^*) = f_1^* \quad \text{and} \quad T(v_n^*) = f_{n+1}^* \text{ for all } n \geq 1,$$

is extended, by Lemma 2 of [1], to a $\sigma(X^*, X)$-continuous isometrical isomorphism from $V$ onto $c^*$.

Let us denote by $^\perp Q$ the annihilator of $Q \subseteq X^*$ in $X$; since $P$ and $T$ are $w^*$-continuous, there exist two linear bounded operators

$$S : {}^{X}\!/_{\perp Q} \to X \quad \text{and} \quad R : c \to {}^{X}\!/_{\perp Q}$$

such that $S^* = P$, $R^* = T$ and $R$ is an isometrical isomorphism. The linear map $S \circ R : c \to X$ is an isometry into. This fact concludes the proof. $\qquad\square$

To conclude this section we point out that, whenever a separable Banach space contains an isometric copy of $c$, then, by Proposition 3.1 in [3], ensures that it also contains a 1-complemented copy of $c$.

## 5. A SPECIAL CLASS OF $\ell_1$-PREDUALS: THE SPACES $W_\alpha$

This section is devoted to study some properties of $W_\alpha$ spaces. We recall that for every $\alpha = (\alpha(1), \alpha(2), \dots) \in B_{\ell_1}$, we define the following hyperplane of $c$:

$$W_\alpha = \left\{ x = (x(1), x(2), \dots) \in c : \lim_{i \to \infty} x(i) = \sum_{i=1}^{\infty} \alpha(i) x(i) \right\}.$$

Then, $W_\alpha^*$ is isometrically identified with $\ell_1$ under the following duality action: for every $x^* \in W_\alpha^*$ there is a unique $f = (f(1), f(2), \dots) \in \ell_1$ such that

$$x^*(x) = \sum_{i=1}^\infty f(i)x(i) = f(x)$$

for every $x = (x(1), x(2), \dots) \in W_\alpha$ (see Theorem 4.3 in [2]). These spaces describe completely, from an isometric point of view, all the $\ell_1$-preduals such that the standard basis of the dual is $w^*$-convergent. Indeed, $X$ is an $\ell_1$-predual such that the standard basis $\{e_n^*\}$ of $X^*$ is $\sigma(X^*, X)$-convergent to $\alpha \in B_{X^*}$ if and only if $X$ is isometric to $W_\alpha$ (see [2]). We also recall when two spaces $W_\alpha$ and $W_\beta$ are isometric.

*Remark* 5.1. Let $\alpha, \beta \in B_{\ell_1}$. Then the spaces $W_\alpha$ and $W_\beta$ are isometric if and only if there exists a finite sequence of signs $\{\epsilon_n\}_{n=1}^{j_0}$, $\epsilon_n = \pm 1$, and a permutation $\pi : \mathbb{N} \to \mathbb{N}$ such that

$$\begin{cases} \alpha(n) = \epsilon_n \beta(\pi(n)) & \text{for } 1 \le n \le j_0 \\ \alpha(n) = \beta(\pi(n)) & \text{for } n > j_0 \end{cases}$$

(see Proposition 2.6 in [33]).

Finally, we characterize $W_\alpha$ spaces that are $A(S)$ spaces of affine continuous functions on a Choquet simplex $S$ (we recall that these spaces are examples of $L_1$-preduals and that were widely studied in the literature; see, e.g., [23] and Section 9 in the present paper).

*Remark* 5.2. Let us consider a space $W_\alpha$, then the following properties are equivalent:

(i) $W_\alpha$ is isometric to a space $A(S)$ of affine continuous functions on a Choquet simplex $S$;
(ii) there exists $\bar{x} \in S_{W_\alpha}$ such that $|\alpha(\bar{x})| = \|\alpha\| = 1$;
(iii) ext $B_{W_\alpha} \ne \emptyset$.

The proof of this equivalences follows from [36], where the author proved that an $L_1$-predual $X$ is an $A(S)$ space if and only if ext $B_X \ne \emptyset$. The details are left to the reader.

The sequel of the section is divided into two subsections. The first one is devoted to the cubic structure of $W_\alpha$ and to finding a sufficient condition to ensure that an $\ell_1$-predual are isometric to a space $W_\alpha$ with $\|\alpha\| = 1$ by means of their cubic structure. This result is worth to be mentioned because the spaces $W_\alpha$ with $\|\alpha\| = 1$ will be relevant in the sequel. The second subsection deals with the conditions ensuring that a space $W_\alpha$ can be isometrically embedded in a space $W_\beta$. These results will play a key role in the study the structure of $\ell_1$-predual (see Section 6).

## 5.1. **Cubic structure of the spaces** $W_\alpha$. The matter contained in this subsection is not strictly useful in the sequel of this paper. However, it seems interesting in itself since it provides an explicit example of the cubic structure, presented in Subsection 5.1 for a generic separable $L_1$-predual.

Let us consider the vectors of a given $W_\alpha$:

$$e_i^n = \left(0, 0, \cdots, 0, \underbrace{1}_{j}, 0, \cdots, 0, \underbrace{a_i^n}_{n+1}, a_i^n, \cdots\right)$$

where, $i, n \in \mathbb{N}$, $1 \leq i \leq n$, $a_i^n = \frac{\alpha(i)}{1 - A_{n+1}}$ and $A_n = \sum_{j=n}^{+\infty} \alpha(j)$. Indeed, it holds for every $1 \leq i \leq n$

$$e_i^n = e_i^{n+1} + a_i^n e_{n+1}^{n+1}.$$

Therefore, the infinite triangular matrix $A = (a_i^n)$ represents the space $W_\alpha$. It is easy to see that the limits of the rows of the matrix $A$ give the components of the element $\alpha$, i.e., we have

$$\lim_{n \to \infty} a_i^n = \alpha(i).$$

Now, we recall from ([30]) a result concerning the evaluation of the functionals $\left\{\phi_j\right\}_{j=1}^\infty$, as defined in (2.4).

*Proposition* 5.3. Let $A = (a_i^n)$ a triangular representation matrix of a separable $L_1$ predual space $X$. Then

(1) $\phi_j(e_k^k) = \Delta_j^{j-k}$ for $j > k$

(2) $\Delta_k^n = a_k^{n+k-1} + \sum_{j=k+1}^{n-1} a_j^{n+k-1} \Delta_k^{j-k}$

where:

$$\Delta_n^k = \det \begin{bmatrix} a_n^n & a_n^{n+1} & a_n^{n+2} & \cdots & a_n^{n+k-1} \\ -1 & a_{n+1}^{n+1} & a_{n+1}^{n+2} & \cdots & a_{n+1}^{n+k-1} \\ 0 & -1 & a_{n+2}^{n+2} & \cdots & a_{n+2}^{n+k-1} \\ \vdots & \vdots & \ddots & \ddots & \vdots \\ 0 & \cdots & \cdots & -1 & a_{n+k-1}^{n+k-1} \end{bmatrix}$$

The previous result allow us to characterize the spaces $W_\alpha$ among the separable $\ell_1$-preduals by using the representing matrices

*Remark* 5.4. By taking in account the definition of the functionals $\left\{\phi_j\right\}_{j=1}^\infty$ given in equation (2.4), it is easily seen that an $\ell_1$-predual $X$ is isometric to a space $W_\alpha$ if and only if

$$\phi_n \xrightarrow[n \to \infty]{\sigma(\ell_1, X)} \sum_{j=1}^{+\infty} \alpha_j \phi_j$$

or, equivalently,

(5.1) $$\lim_{n \to \infty} \phi_n(e_k^k) = \sum_{j=1}^{+\infty} \alpha(j)\phi_j(e_k^k) = \alpha_k + \sum_{j=k+1}^{+\infty} \alpha(j)\phi_j(e_k^k).$$

By Proposition 5.3, an equivalent formulation of (5.1) is

(5.2) $$\lim_{n \to \infty} \Delta_k^n = \alpha(k) + \sum_{j=k+1}^{+\infty} \alpha(j)\Delta_k^{j-k} \text{ for } k = 1, 2, \cdots.$$

The last equation provides a necessary and sufficient condition to guarantee that an $\ell_1$-predual $X$ is isometric to $W_\alpha$ in terms of its representing matrix.

Now, we focus our attention to the spaces $W_\alpha$ with $\|\alpha\| = 1$. Indeed, we prove a sufficient condition to ensure that an $\ell_1$-preduals is isometric to the space $W_\alpha$, when $\|\alpha\| = 1$. Suppose that $X^* = \ell_1$ and let $A = (a_{j,n})$ be a triangular matrix which represents the space $X$ with $\lim_{n\to\infty} a_{j,n} = \alpha_j$. If $\sum_{j=1}^{+\infty} |\alpha_j| = 1$, then $X$ is isometric to $W_\alpha$.

*Proposition 5.5.* Suppose that $X$ is an $\ell_1$-predual space and let $A = (a_i^n)$ be a triangular matrix which represents the space $X$ such that $\lim_{n\to\infty} a_i^n = \alpha(i)$. If $\|\alpha\| = \sum_{j=1}^{+\infty} |\alpha(j)| = 1$, then $X$ is isometric to $W_\alpha$.

*Proof.* By Remark 5.4, it suffices to show that

$$(5.3) \qquad \lim_{n\to\infty} \Delta_k^n = \alpha(k) + \sum_{j=k+1}^{+\infty} \alpha(j)\Delta_k^{j-k} \text{ for } k = 1, 2, \cdots.$$

Choose $\epsilon > 0$ and let $N$ be such that $\sum_{j=1}^n |\alpha(j)| \geq 1 - \epsilon$ for every $n \geq N$ and let $M > N$ be such that for every $n > M$ and $j = 1, 2, \ldots, N$ we have $|a_j^n - \alpha(j)| \leq \epsilon/N$. For such $n$ we have:

$$1 \geq \sum_{j=1}^n |a_j^n| \geq \sum_{j=1}^N |\alpha(j)| - \sum_{j=1}^N |\alpha(j) - a_j^n| + \sum_{j=N+1}^n |a_j^n| \geq 1 - 2\epsilon + \sum_{j=N+1}^n |a_j^n|$$

From this inequality it follows that $\sum_{j=N+1}^n |a_j^n| \leq 2\epsilon$ for every $n > M$. Again let $n > M$, then:

$$(5.4) \qquad \sum_{j=1}^n |\alpha(j) - a_j^n| \leq \sum_{j=1}^N |\alpha(j) - a_j^n| + \sum_{j=N+1}^n |a_j^n| + \sum_{j=N+1}^n |\alpha(j)| \leq 4\epsilon$$

Now, let us go back to (5.3). Fix $k \geq 1$ and $n > M$ such that $|\alpha(k) - a_k^{n+k-1}| < \epsilon$. Then,

$$\left| \Delta_k^n - \left( \alpha(k) + \sum_{j=k+1}^{+\infty} \alpha(j)\Delta_k^{j-k} \right) \right| = \left\| \left( a_k^{n+k-1} + \sum_{j=k+1}^{n-1} a_j^{n+k-1}\Delta_k^{j-k} \right) - \left( \alpha(k) + \sum_{j=k+1}^{+\infty} \alpha(j)\Delta_k^{j-k} \right) \right\| \leq$$

$$\leq |\alpha_k - a_k^{n+k-1}| + \sum_{j=k+1}^{n-1} \left| \left( \alpha(j) - a_j^{n+k-1} \right) \Delta_k^{j-k} \right| + \sum_{j=n}^{+\infty} |\alpha(j)\Delta_k^{j-k}| \leq$$

$$\leq \epsilon + \sum_{j=1}^{n+k-1} \left| \alpha(j) - a_j^{n+k-1} \right| + \sum_{j=n}^{+\infty} |\alpha(j)| \leq \epsilon + 4\epsilon + \epsilon = 6\epsilon.$$

□

### 5.2. Isometric embeddings of $W_\alpha$ in $W_\beta$.
The present subsection is devoted to answer the question: *given a space $W_\alpha$, what are the elements $\beta \in B_{\ell_1}$ such that $W_\alpha$ contains a subspace isometric to $W_\beta$?*

Our first result shows that it is always possible to embed into a generic $W_\alpha$ an isometric copy of $W_\beta$ with $\beta = (r, 0, 0, \ldots)$ for every $r < \|\alpha\|$. We point out that these spaces are $M$-spaces (see Section 9) for every $r > 0$ and they will play a key role in our approach (see Section 6 and 7).

*Proposition* 5.6. Let $\alpha \in B_{\ell_1}$ be such that $\|\alpha\| > 0$. Then, $W_\alpha$ contains an isometric copy of $W_{re_1^*}$, for every $r \in \mathbb{R}$ such that $|\alpha(1)| \leq r < \|\alpha\|$.

*Proof.* By taking in account Remark 5.1, without loss of generality, we can consider only the case where $\alpha$ is such that

$$\alpha(1) \geq |\alpha(n+1)| \text{ for every } n \geq 1.$$

From the assumption about $r$, it is easy to see that there exists an integer number $K$, $K > 1$ such that

$$\sum_{j=1}^{K-1} |\alpha(j)| \leq r < \sum_{j=1}^{K} |\alpha(j)|.$$

Therefore, there exists a real number $s \in [0, 1)$ such that $r = s|\alpha(1)| + \sum_{j=2}^{K} |\alpha(j)|$.

Now, let $N > K$ be an integer number such that

$$(5.5) \qquad \sum_{j=N}^{\infty} |\alpha(j)| \leq (1-s)|\alpha(1)|$$

and let $w = (w(1), w(2), \dots) \in W_{re_1^*}$, then we define the map $T : W_{re_1^*} \to W_\alpha$ by

$$T(w)(i) = \begin{cases} sw(1) - \sum_{j=N}^{\infty} \frac{\alpha(j)}{\alpha(1)} w(j-N+2) & \text{if} & i = 1 \\ \text{sgn}(\alpha(i))w(1) & \text{if} & 2 \leq i \leq K \\ 0 & \text{if} & K+1 \leq i \leq N-1 \\ w(i-N+2) & \text{if} & i \geq N. \end{cases}$$

By an easy computation we obtain that

$$(5.6) \qquad \sum_{j=1}^{\infty} T(w)(j) = rw(1).$$

On the other hand, since $w \in W_{re_1^*}$, by definition of the map $T$ we have that

$$(5.7) \qquad \lim_n T(w)(n) = \lim_n w(n) = rw(1).$$

Therefore, by (5.6) and (5.7), we have $\sum_{j=1}^{\infty} T(w)(j) = \lim_n T(w)(n)$, hence $T(w) \in W_\alpha$. If we show that $T$ is an isometry we have proved our proposition. First, since the norms in $W_\alpha$ and in $W_{re_1^*}$ are the sup norm, immediately we have that $\|T(w)\| \geq \|w\|$ for every $w \in W_{re_1^*}$. Second, by recalling (5.5), we obtain

$$\left| sw(1) - \sum_{j=N}^{\infty} \frac{\alpha(j)}{\alpha(1)} w(j-N+2) \right| \leq \left( s + \frac{\sum_{j=N}^{\infty} |\alpha(j)|}{|\alpha(1)|} \right) \|w\| \leq \|w\|.$$

Hence, it holds $\|T(w)\| \leq \|w\|$ for every $w \in W_{re_1^*}$. $\qquad \square$

Our next results show that Proposition 5.6 is no longer true even if we take $r = \|\alpha\|$. The case $r = \|\alpha\| = 1$ follows from Proposition 2.1 in [3]:

*Proposition* 5.7. Let $\alpha \in \ell_1$ be such that $\|\alpha\| = 1$. Then the following are equivalent.

(1) $W_\alpha$ contains an isometric copy of $W_{e_1^*}$.

(2) The set $\{n \in \mathbb{N} : \alpha(n) < 0\}$ is finite and $\{n \in \mathbb{N} : \alpha(n) = 0\}$ is infinite.

If $r = \|\alpha\| < 1$, then we have the following characterization.

*Proposition* 5.8. Let $\alpha \in \ell_1$ be such that $0 < \|\alpha\| < 1$ and let $r = \|\alpha\|$. Then the following are equivalent.

(1) $W_\alpha$ contains an isometric copy of $W_{re_1^*}$.

(2) supp $(\alpha)$ is finite.

*Proof.* (2) $\implies$ (1). Suppose that supp $(\alpha)$ is finite. By Remark 5.1, we can assume that supp $(\alpha) = \{1, \dots, n\}$ for some $n \in \mathbb{N}$. Let $T : W_{re_1^*} \to W_\alpha$ be defined for every $x = (x(1), x(2), \dots) \in W_{re_1^*}$ by

$$(Tx)(i) = \begin{cases} \mathrm{sgn}(\alpha(i))x(1) & \text{for } i \in \mathrm{supp}\,(\alpha) \\ x(i-n) & \text{for } i \geq n+1. \end{cases}$$

It is easy to check that $T$ is an isometry into.

(1) $\implies$ (2). Suppose that $W_\alpha$ contains an isometric copy of $W_{re_1^*}$. We have

$$e_i^* \xrightarrow[i \to \infty]{\sigma\left(\ell_1, W_{re_1^*}\right)} re_1^*.$$

By Lemma 3.1, every $e_i^* \in W_{re_1^*}^*$ admits an extension $\widetilde{e_{k_i}^*}$ to the whole space $W_\alpha$ such that $\widetilde{e_{k_i}^*} \in \mathrm{ext}\, B_{W_\alpha^*}$. W.L.O.G. we can assume that

$$\widetilde{e_{k_i}^*} \xrightarrow[i \to \infty]{\sigma(\ell_1, W_\alpha)} \alpha.$$

Then $\alpha\big|_{W_{re_1^*}} = re_1^*$. Observe now that there exists $\bar{x} \in S_{W_{re_1^*}}$ such that

$$re_1^*(\bar{x}) = \big\|re_1^*\big\| = r,$$

hence, $\alpha(\bar{x}) = r = \|\alpha\|$. This implies that supp $(\alpha)$ must be a finite set. Indeed, if supp $(\alpha)$ is infinite and $\|\alpha\| < 1$, then $\alpha \in W_\alpha^*$ does not attain its norm. $\qquad\square$

Now we are in a position to present the following

*Theorem* 5.9. For every $\alpha \in B_{\ell_1}$, $W_\alpha$ contains an isometric copy of $W_\beta$ whenever $\beta \in \ell_1$ and $\|\beta\| < \|\alpha\|$.

*Proof.* By Proposition 5.6, $W_\alpha$ contains an isometric copy of $W_{re_1^*}$ with $\|\beta\| \leq r < \|\alpha\|$. Let us consider the map $T : W_\beta \to W_{re_1^*}$ defined by

$$T(x) = \left(\frac{x_\infty}{r}, x(1), x(2), \dots\right),$$

where $x = (x(1), x(2), \dots) \in W_\beta$ and $x_\infty = \lim_n x(n) = \sum_{i=1}^\infty \beta(i)x(i)$. It is easy to see that $T$ is an injective isometry. $\qquad\square$

It is also worth noting that, if we choose arbitrarily $\alpha$ and $\beta$ with the same norm, in general neither $W_\beta$ is a subspace of $W_\alpha$ nor $W_\alpha$ of $W_\beta$. For example, let

$$\alpha = (1/2, 1/4, 1/8, 1/16, \dots)$$

and

$$\beta = (1/2, -1/4, 1/8, -1/16, \dots).$$

First, we observe that $W_\beta$ does not contain an isometric copy of $W_\alpha$. Indeed, $W_\alpha$ is an $A(S)$ space (see Section 9), whereas $W_\beta$ does not contain any $A(S)$ space (see Example 3.3 in [4]). Moreover, $W_\alpha$ does not contain an isometric copy of $W_\beta$, which is a simple consequence of the following result.

*Proposition* 5.10. Let $\alpha \in S_{\ell_1}$ with $\alpha(n) > 0$ for every $n \in \mathbb{N}$ and suppose that $W_\alpha$ contains an isometric copy of $W_\beta$ with $\beta \in S_{\ell_1}$. Then $\{n \in \mathbb{N} : \beta(n) < 0\}$ is finite and $\{n \in \mathbb{N} : \beta(n) = 0\}$ is empty.

*Proof.* Let $(e_n^*)$ be the standard basis in $W_\beta^*$ and suppose that

$$A = \{n \in \mathbb{N} : \beta(n) < 0\}$$

is infinite. By Lemma 3.1, let $f_n^* \in \operatorname{ext} B_{W_\alpha^*}$ be such that $f_n^*\big|_{W_\beta} = e_n^*$, for every $n \in A$. W.L.O.G. we may assume that

$$f_n^* \xrightarrow{\sigma(\ell_1, W_\alpha)} \alpha.$$

Then, $\alpha\big|_{W_\beta} = \beta$. Therefore, for sufficiently large $n \in A$, we have

$$2 = \|\beta - e_n^*\| \leq \|\alpha - f_n^*\| < 2.$$

A contradiction. Suppose now that

$$B = \{n \in \mathbb{N} : \beta(n) = 0\}$$

is nonempty and let $n_0 \in B$. By Lemma 3.1, there exists $f^* \in \operatorname{ext} B_{W_\alpha^*}$ such that $f^*\big|_{W_\beta} = e_{n_0}^*$. Then, $-f^*\big|_{W_\beta} = -e_{n_0}^*$. By following the reasoning from the first part of the proof, we conclude that there exists $g^* \in \{\pm\alpha\}$ such that $g^*\big|_{W_\beta} = \beta$. Then, $\|\pm e_{n_0}^* - \beta\| = 2$, but $\|f^* - g^*\| < 2$ or $\|-f^* - g^*\| < 2$. Again a contradiction. $\quad\square$

On the other hand, it may also happen that $W_\alpha$ and $W_\beta$ are both isometrically embeddable one in the other, although none of them contains an isometric copy of $c$, as shown in the following example.

*Example* 5.11. Let $\alpha = \left(\frac{1}{n(n+1)}\right)$ and $\beta = \left(\frac{1}{2^n}\right)$. Define

$$Sw = \left(\sum_{n=1}^{+\infty} \gamma_n w(n), w(1), w(2), \dots\right),$$

where $w = (w(1), w(2), \dots) \in W_\alpha$ and $\gamma_n = 2\,(\alpha(n) - \beta(n+1))$. Then $S : W_\alpha \to W_\beta$ is an isometric embedding of $W_\alpha$ into $W_\beta$. Moreover, we can also embeds isometrically $W_\beta$ into $W_\alpha$ via the following isometry:

$$Tw = (\overbrace{w(1)}^{1}, \overbrace{w(2), w(2)}^{2}, \overbrace{w(3), w(3), w(3), w(3)}^{4}, \dots) \quad \text{for } w \in W_\beta.$$

Finally, let us consider the norm one projection $P : W_\alpha^* \to W_\alpha^*$ defined by

$$Px^* = \left(x^*(1), \underbrace{\frac{1}{2}\sum_{i=2}^{3} x^*(i), \frac{1}{2}\sum_{i=2}^{3} x^*(i)}_{2}, \underbrace{\frac{1}{4}\sum_{i=4}^{7} x^*(i), \dots, \frac{1}{4}\sum_{i=4}^{7} x^*(i)}_{4}, \dots\right).$$

Then, we have $P(x^*)(Tw) = x^*(Tw)$ for every $w \in W_\beta$.

The last part of this subsection is devoted to a discussion of some properties about embeddings among $W_\alpha$ spaces, that point out connections of our approach with some old results. Indeed, the existence of the projection explicitely stated at the end of Example 5.11 is guaranteed by a more general result (Lemma 3, §22 in [23]).

*Proposition* 5.12. Let $X$ and $Y$ two $L_1$-preduals such that $Y \subset X$. Then there exists a norm one projection $P : X^* \to X^*$ whose kernel is $Y^\perp$ and whose range is isometric to $Y^*$ under the restriction mapping $x^* \to x^*|_Y$ for $x^* \in X^*$.

It is well known that the norm one projection from $\ell_1$ to $\ell_1$ has a concrete representation (see, e.g., [28]).

*Proposition* 5.13. Let $P : \ell_1 \to \ell_1$ be a norm one projection. Then, there exist a collection of subsets $\{\sigma_j\}_{j\in\mathbb{N}}$ of $\mathbb{N}$ such that $\sigma_j \cap \sigma_k = \emptyset$ for $j \neq k$ and norm one vectors $\{u_j^*\}_{j\in\mathbb{N}}$ of the form: $u_j^* = \sum_{i\in\sigma_j} \lambda_i e_i^*$ so that $Px^* = \sum_{j=1}^{\infty} u_j^{**}(x^*) u_j^*$, where $\{u_j^{**}\}_{j\in\mathbb{N}} \subset \ell_\infty$ and the vectors $u_j^{**}$ satisfy $\|u_j^{**}\| = u_j^{**}(u_j^*) = 1$ for every $j \in \mathbb{N}$.

Some other properties, whose proofs are straightforward, can be added:

(a) $u_j^{**}(u_k^*) = 0$ for every $j, k \in \mathbb{N}$ such that $j \neq k$.
(b) If $k \in \sigma_j$ then $u_j^{**}(e_k^*) = \text{sgn}(\lambda_k)$.
(c) Let $\omega_j = \text{supp}(u_j^{**})$. Then $\sigma_j \subset \omega_j$ and if $k \neq j$ then $\omega_k \cap \sigma_j = \emptyset$.
(d) $P(x^*) = 0$ if and only if $u_j^{**}(x^*) = 0$ for every $j \in \mathbb{N}$

Now, let us suppose that there exists an isometric embedding $T$ of $W_\beta$ in $W_\alpha$ and, according to the previous propositions, let $P$ be the norm one projection from $W_\alpha^* = \ell_1$ on $W_\alpha^*$. In particular, if $\{e_n^*\}$ denote the canonical basis of $W_\alpha^*$, then for every $s \in \sigma_j$ we have

$$(Tw)(s) = e_s^*(Tw) = P(e_s^*)(Tw) = \text{sgn}(\lambda_s) \sum_{i\in\sigma_j} \lambda_i (Tw)(i).$$

Moreover, Lemma 3.1 implies that, if $\{f_n^*\}$ is the canonical basis of $W_\beta^*$, then there exists $k_n \in \mathbb{N}$ such that for every $w \in W_\beta$:

$$(Tw)(k_n) = e_{k_n}^*(Tw) = \pm f_n^*(w) = \pm w(n).$$

Now, we assume that $k_n \in \sigma_{j_n}$. By the previous facts, it follows that for every $s \in \sigma_{j_n}$:

$$(5.8) \qquad\qquad (Tw)(s) = \pm w(n)$$

where, of course, the sign of $w(n)$ depends on $s$. Moreover, since $T$ is a linear operator, for $s \in \mathbb{N}$ we have

$$(Tw)(s) = \sum_{n=1}^{\infty} \gamma_n w(n) = \sum_{n=1}^{\infty} \pm \gamma_n (Tw)(k_n) = \left( \sum_{n=1}^{\infty} \pm \gamma_n e_{k_n}^* \right)(Tw).$$

The last equality implies that $P^*\left(e_s^* - \sum_{n=1}^{\infty} \pm \gamma_n e_{k_n}^*\right) = 0$ or, equivalently, for every $j \in \mathbb{N}$,

$$u_j^{**}\left(e_s^* - \sum_{n=1}^{\infty} \pm \gamma_n e_{k_n}^*\right) = 0.$$

By recalling properties (c) above, we obtain that

$$u_{j_n}^{**}(e_s^*) = \pm\gamma_n,$$

hence that $s \in \omega_{j_n}$ for every $n$ such that $\gamma_n \neq 0$.

Summing up, if we have the projection $P : W_\alpha^* \to W_\alpha^*$, from Proposition 5.13 and the relationship, where $k_n$ corresponds to $n$, between $\{f_n^*\}$ and $\{e_{k_n}^*\}$ given by Lemma 3.1, we get the description of the embedding $T$ and we can also obtain the value of $\beta$ evaluating:

$$\beta(k) = \sum_{h \in \omega_{j_k}} u_{j_k}^{**}(e_h^*)\alpha(h).$$

*Example* 5.14. (Continuation of Example 5.11) Let $\alpha_n = \frac{1}{2^n}$. If $P : W_\alpha^* \to W_\alpha^*$ be such that $\sigma_j = \{j + 1\}$, $\omega_k = \{1, k + 1\}$,

$$u_j^{**} = \left( \left( \frac{2}{j(j+1)} - \frac{1}{2^j} \right), 0, \ldots, \underbrace{1}_{j+1}, 0, \ldots \right),$$

and $k_n = n$, then we obtain $\beta_n = \frac{1}{n(n+1)}$ and the isometry $S$ of Example 5.11.

Finally, we point out also that formula (5.8) shows that, if there is an isometric copy of $W_\beta$ inside $W_\alpha$, there always exists an isometric embedding $T : W_\beta \to W_\alpha$ such then each component $w(n)$ of $w \in W_\beta$ must appear, up to the sign, as a component of $Tw \in W_\alpha$.

## 6. Subspaces of $\ell_1$-preduals

This section is the core of our paper. Here we show that the features of the $\sigma(\ell_1, X)$-topology provide information on the subspaces of $X$ in terms of the presence of suitable $W_\alpha$ spaces.

In [5], we proved a result concerning the relationship between $w^*$-cluster points of the standard basis of $\ell_1$ and the presence of a $W_\alpha$ space as a subspace of $X$ (and something more). We point out that the assumptions of this proposition require the knowledge of the behaviour of the whole subsequence of the standard basis of $\ell_1$ converging to a $w^*$-cluster point.

*Theorem* 6.1. (Proposition 3.1 in [5]) Let $X$ be an $\ell_1$-predual. If there exists a subsequence $\{e_{n_k}^*\}_{k=1}^\infty$ and an element $\alpha \in B_{\ell_1}$ such that

$$e_{n_k}^* \xrightarrow[k \to \infty]{\sigma(\ell_1, X)} \alpha,$$

and $(\operatorname{supp}(\alpha) \setminus \{n_k\}_{k=1}^\infty) = \{i_1, \ldots, i_m\}$ is a finite set, then $X$ contains a 1-complemented isometric copy of $W_\beta$ with

$$\beta = (\alpha(i_1), \ldots, \alpha(i_m), \alpha(n_1), \alpha(n_2), \ldots, \alpha(n_k), \ldots).$$

A natural question arises whether we may conclude the presence of $W_\alpha$, knowing only the set of $w^*$-cluster points of the standard basis in $\ell_1$. Actually, the following theorem shows that information about radius of this set is necessary and sufficient to determine the existence of subspaces isometric to certain $W_\alpha$ spaces. Therefore, this theorem can be also seen as a "quantitative" version of Theorem 1 in [38]. Its proof follows the lines of the proof of Theorem 1 in [38], by introducing some modifications in order to specify the type of $W_\alpha$ space contained in $X$.

*Theorem* 6.2. For every $\ell_1$-predual $X$ we have:

$$r^*(X) = \sup\left\{\|g^*\| : g^* \in (\text{ext } B_{\ell_1})'\right\} = \sup\left\{\|\alpha\| : \alpha \in B_{\ell_1}, W_\alpha \subset X\right\}.$$

*Proof.* Without loss of generality, we may restrict our attention to the set of $w^*$-cluster points of the standard basis in $\ell_1$.

We begin by proving that $r^*(X) \leq \sup\{\|\alpha\| : \alpha \in B_{\ell_1}, W_\alpha \subset X\}$. Observe first that for every $\varepsilon > 0$ there exist $g^* \in (\text{ext } B_{\ell_1})'$ and $\bar{x}_\varepsilon \in X$ such that

(1) $\|\bar{x}_\varepsilon\| = 1$,
(2) $|g^*(\bar{x}_\varepsilon)| \geq (1 - \varepsilon)r^*(X)$.

We consider the "cubic" structure $\{X_n\}_{n=1}^\infty$ of $X$ as outlined in Subsection 2.1. By Remark 2.2, we can choose $X_1 = \text{span}(\bar{x}_\varepsilon)$. Now, let us consider the functionals $\left\{\phi_j\right\}_{j=1}^\infty$ as defined in (2.4). Lemma 2.1 implies that there exists a subsequence $\left\{\phi_{j_k}\right\}_{k=1}^\infty$ such that

$$(6.1) \qquad\qquad \phi_{j_k} \xrightarrow[k\to\infty]{\sigma(X^*,X)} g^*.$$

Without loss of generality we can assume that $j_1 = 1$, hence $\phi_{j_1}(\bar{x}_\varepsilon) = 1$.

We are now in a position to use Lemma 2.3 to build up a subspace $U$ of $X$ such that, for every $u \in U$

$$(6.2) \qquad\qquad \lim_{s\to\infty} \psi_s(u) = g^*(u),$$

where $\psi_s = \phi_{j_s}|_U$ for every $s = 1, , 2, \ldots$. Let $\psi_0 = g^*|_U$.

The proof falls naturally into two cases. First we consider $\|g^*\| > 0$, then we will move to the case $\|g^*\| = 0$.

If $\|g^*\| > 0$ then $\psi_0(\bar{x}_\varepsilon) \neq 0$, hence $\psi_0 \neq 0$. Moreover, by (2.8), we have

$$(6.3) \qquad\qquad \|u\| = \sup_{s\geq 0} |\psi_s(u)| \quad \text{for every } u \in U.$$

It is worth noting that, if $u = \bar{x}_\varepsilon$, the previous inequality gives

$$(6.4) \qquad\qquad 1 = \|\bar{x}_\varepsilon\| = \sup_{s\geq 0} |\psi_s(\bar{x}_\varepsilon)| = \psi_1(\bar{x}_\varepsilon).$$

Relation (6.3) makes obvious that the linear map $T : U \to c$ defined by $T(u) = (\psi_1(u), \psi_2(u), \psi_3(u), \dots)$ is an isometry from $U$ onto a subspace $Y$ of $c$. If $T$ is onto $c$, we get the conclusion. Now, let us suppose that $T$ is not onto $c$. Then, there is a non null functional $e^* \in c^*$ such that $T(U) \subseteq \ker(e^*)$. By recalling the notations for $c^*$ introduced in the Introduction, we have a representation $e^* = \sum_{i=1}^\infty a_i e_i^*$. Moreover, since for every $u \in U$

$$T^*\left(e_j^*\right)(u) = e_j^*(T(u)) = e_j^*(\psi_1(u), \psi_2(u), \psi_3(u), \dots) = \psi_{j-1}(u),$$

we have that $T^*\left(e_j^*\right) = \psi_{j-1}$ for every $j \geq 1$ and, in consequence, $T^*(e^*) = \sum_{j=1}^\infty a_j \psi_{j-1}$.

The following equalities hold:

$$(6.5) \qquad 0 = T^* (e^*) (u) = \sum_{i=1}^{\infty} a_i \psi_{i-1} (u) = a_1 \psi_0 (u) + \sum_{i=2}^{\infty} a_i \psi_{i-1} (u)$$

for every $u \in U$. By Lemma 2.3, the equality $\psi_{i-1} \left( u_j^m \right) = \delta_{i-1,j}$ holds for every $m \geq 1$ and $2 \leq i \leq m+1$. Consequently, if $u = u_j^m$, (6.5) yields:

$$(6.6) \quad 0 = T^* (e^*) \left( u_j^m \right) = \sum_{i=1}^{\infty} a_i \psi_{i-1} \left( u_j^m \right) = a_1 \psi_0 \left( u_j^m \right) + a_{j+1} + \sum_{i=m+2}^{\infty} a_i \psi_{i-1} \left( u_j^m \right).$$

Since $e^* \in c^* = \ell_1$, the series $\sum_{i=1}^{\infty} |a_i|$ converges. Therefore, we get from (6.6) that

$$(6.7) \qquad \left| a_1 \psi_0 \left( u_j^m \right) + a_{j+1} \right| = \left| \sum_{i=m+2}^{\infty} a_i \psi_{i-1} \left( u_j^m \right) \right| \xrightarrow[m \to \infty]{} 0.$$

We observe that $a_1 \neq 0$. Indeed, If $a_1 = 0$, from (6.7) we conclude that $a_j = 0$ for every $j \geq 2$, hence that $e^* = 0$, and finally that $T$ is onto. Since $e^* \neq 0$, there is no loss of generality in assuming $\|e^*\| = 1$. First we show that the functional $e^*$ is uniquely defined, up to a numerical multiple (to do this we proceed analogously to the proof of Theorem 1 in [38]), second we prove that $X$ contains an isometric copy of a suitable $W_\beta$. Relation 6.7 implies that $\lim \psi_0 \left( u_j^m \right)$ exists and

$$\frac{a_{j+1}}{a_1} = - \lim_{m \to \infty} \psi_0 \left( u_j^m \right),$$

for all $j \geq 1$, which proves that all $e^*$ with $a_1 \neq 0$ such that $T(U) \subseteq \ker(e^*)$ are numerical multiples of each other, hence $T(U)$ is a hyperplane of $c$. Now, we recall that $u_1^1 = \bar{x}_\varepsilon$, hence (6.6) gives

$$(6.8) \qquad 0 = T^* (e^*) (\bar{x}_\varepsilon) = \sum_{i=1}^{\infty} a_i \psi_{i-1} (\bar{x}_\varepsilon) = a_1 \psi_0 (\bar{x}_\varepsilon) + a_2 + \sum_{i=3}^{\infty} a_i \psi_{i-1} (\bar{x}_\varepsilon).$$

As $a_1 \neq 0$, $\psi_0 = g^*|_U$ and $\bar{x}_\varepsilon \in U$, the previous relation implies that

$$(6.9) \qquad g^* (\bar{x}_\varepsilon) = -\frac{1}{a_1} \left( a_2 + \sum_{i=3}^{\infty} a_i \psi_{i-1} (\bar{x}_\varepsilon) \right).$$

and hence, by the definition of $\bar{x}_\varepsilon$, we get

$$(6.10) \qquad (1 - \varepsilon) \|g^*\| \leq |g^*(\bar{x}_\varepsilon)| \leq \sum_{i=2}^{\infty} \left| \frac{a_i}{a_1} \right|.$$

The hyperplane $T(U)$ of $c$ can be seen as
$$(6.11)$$
$$T(U) = \ker (e^*) = \left\{ w = (w(1), w(2), \dots) \in c : \lim_{i \to \infty} w(i) = \sum_{i=2}^{\infty} \left( -\frac{a_i}{a_1} w(i-1) \right) \right\}.$$

Moreover, according to (6.5) and by recalling that $\left\{\psi_j\right\}_{j=1}^\infty$ is isometrically equivalent to the standard basis of $\ell_1$, we have

$$\|\psi_0\| = \left\| -\frac{1}{a_1} \sum_{i=2}^\infty a_i \psi_{i-1} \right\| = \sum_{i=2}^\infty \left| \frac{a_i}{a_1} \right| \le \|g^*\|,$$

last inequality being a consequence of our definition $\psi_0 = g^*|_U$. If we write $\beta_\varepsilon = (-a_2/a_1, -a_3/a_1, \dots) \in \ell_1$, we see at once that

$$(1 - \varepsilon)\|g^*\| \le \|\beta_\varepsilon\| \le \|g^*\| \le 1$$

and $T(U) = W_{\beta_\varepsilon}$.

The proof is completed by considering the case where $\|g^*\| = 0$. Then, the isometry $T$ defined as above is an isometry onto $c_0$, than can be seen as $W_\beta$ with $\beta = 0$.

Now we show that $r^*(X) \ge \sup\{\|\alpha\| : \alpha \in B_{\ell_1}, \ W_\alpha \subset X\}$. Suppose that there exists $\alpha \in B_{\ell_1}$ with $\|\alpha\| > r^*(X)$ such that $W_\alpha \subset X$. We have $e_n^* \xrightarrow[n\to\infty]{\sigma(\ell_1, W_\alpha)} \alpha$. By applying Lemma 3.1, we conclude that there exists $g^* \in (\mathrm{ext}\, B_{X^*})'$ such that $g^*|_{W_\alpha} = \alpha$. Hence, $\|g^*\| \ge \|\alpha\| > r^*(X)$, a contradiction.  □

Theorem 6.2 has several consequences, that we will develop in the sequel of the paper. We start by the following corollary that can be proved by combining Theorem 6.2 with Theorem 5.9.

*Corollary* 6.3. Let $X$ be an $\ell_1$-predual. Then, $X$ contains an isometric copy of $W_\beta$ for all $\beta \in B_{\ell_1}$ such that $\|\beta\| < r^*(X)$.

*Remark* 6.4. The strict inequality "$<$" in Corollary 6.3 cannot be replaced by "$\le$" even if there is $\alpha \in (\mathrm{ext}\, B_{\ell_1})'$ such that $\|\alpha\| = r^*(X)$. Indeed, the $\ell_1$-predual

$$W_{x_1^*, x_2^*} = \left\{ f \in C(\Omega^2) : f(1,0) = \sum_{j=1}^\infty \frac{f\left(1, \frac{1}{j}\right)}{2^{2j-1}} + \sum_{j=1}^\infty \frac{f\left(\frac{1}{2}, \frac{1}{j}\right)}{2^{2j-2}} \text{ and } f\left(\frac{1}{2}, 0\right) = 0 \right\}$$

is such that

$$e_{1,j}^* \xrightarrow{\sigma(\ell_1(\Omega_0^2), W_{x_1^*, x_2^*})} x_1^* = \left( \frac{1}{2}, \frac{1}{8}, \frac{1}{32}, \dots, \frac{1}{4}, \frac{1}{16}, \frac{1}{64}, \dots \right)$$

and

$$e_{2,j}^* \xrightarrow{\sigma(\ell_1(\Omega_0^2), W_{x_1^*, x_2^*})} x_2^* = (0, 0, 0, \dots, 0, 0, 0, \dots)$$

(for the notations see Example 4.2). It can be proved that there is no $\beta \in \ell_1$ with $\|\beta\| = 1$ such that $W_\beta \subset W_{x_1^*, x_2^*}$ (for the details, see Example 3.4 in [5]).

Now, we pass to consider a sufficient condition to ensure the presence of an isometric copy of the space $W_\alpha$.

*Corollary* 6.5. Let $X$ be an $\ell_1$-predual. Let $g^* \in (\mathrm{ext}\, B_{\ell_1})'$ and suppose that there exists $\bar{x} \in S_X$ such that

$$|g^*(\bar{x})| = \|g^*\|.$$

Then, there exists $\alpha \in B_{\ell_1}$ such that

(i) $\|\alpha\| = \|g^*\|$,

(ii) $X$ contains an isometric copy of $W_\alpha$.

*Proof.* We only need to consider the case when $\|g^*\| > 0$. Following the same notations as in the proof of Theorem 6.2, we point out that $\bar{x} \in U$, hence, equation (6.9) gives

$$\|g^*\| = |g^*(\bar{x})| = \sum_{i=2}^{\infty} \left| \frac{a_i}{a_1} \right|.$$

Since $\alpha = (-a_2/a_1, -a_3/a_1, \dots)$, then $\|\alpha\| = \|g^*\|$, and the proof is completed. □

*Remark* 6.6. Under the assumptions of Corollary 6.5 and by using the same notations as in the proof of Theorem 6.2, it is worth noting that $T(\bar{x}) \in T(U) = W_\beta$. Hence, by using the duality between $W_\beta$ and $\ell_1$ (see Section 5) we get

$$|\beta(T(\bar{x}))| = \left| \sum_{i=2}^{\infty} -\frac{a_i}{a_1} T(\bar{x})(i) \right| = \left| \lim_{i \to \infty} \psi_i(\bar{x}) \right| = |\psi_0(\bar{x})| = |\alpha|_U(\bar{x})| = \|\alpha\| = \|\beta\|.$$

We point out that the existence of an isometric copy of a $W_\alpha$ space in $\ell_1$-predual $X$ was already presented in Theorem 6.1. Both, in Corollary 6.5 and Theorem 6.1, we assume that there exists a $\sigma(X^*, X)$-cluster point $\alpha \in B_{\ell_1}$ of the standard basis of $X^* = \ell_1$, but in Theorem 6.1 we require essentially that:

$$\text{if } e^*_{n_k} \xrightarrow[k \to \infty]{\sigma(\ell_1, X)} \alpha \text{ then supp}(\alpha) \subseteq \{n_k\}_{k=1}^{\infty}, \quad \text{(S)}$$

whereas in Corollary 6.5 the key assumption was

$$\text{there exists } \bar{x} \in S_X \text{ such that } |\alpha(\bar{x})| = \|\alpha\|. \quad \text{(P)}$$

The following two examples prove that properties (P) and (S) are completely unrelated. In this examples, we adopt the notations introduced in Example 4.2.

*Example* 6.7. Let $x_1^*, x_2^* \in \ell_1(\Omega_0^2)$ be defined by

$$x_1^*(i, j) = \begin{cases} 0 & \text{for } i = 1 \text{ and } j \geq 1 \\ \frac{1}{2^j} & \text{for } i = 2 \text{ and } j \geq 1 \end{cases},$$

and

$$x_2^*(i, j) = \begin{cases} \frac{1}{2^j} & \text{for } i = 1 \text{ and } j \geq 1 \\ 0 & \text{for } i = 2 \text{ and } j \geq 1 \end{cases}.$$

Then

$$W_{x_1^*, x_2^*} = \left\{ f \in C(\Omega^2) : f(1, 0) = \sum_{j=1}^{\infty} \frac{f\left(\frac{1}{2}, \frac{1}{j}\right)}{2^j} \text{ and } f\left(\frac{1}{2}, 0\right) = \sum_{j=1}^{\infty} \frac{f\left(1, \frac{1}{j}\right)}{2^j} \right\}.$$

By Theorem 2.1 in [5], $W_{x_1^*, x_2^*}^* = \ell_1(\Omega_0^2) = \ell_1$ and

$$e^*_{1,j} \xrightarrow{\sigma(\ell_1(\Omega_0^2), W_{x_1^*, x_2^*})} x_1^* = \left(0, 0, 0, \dots, \frac{1}{2}, \frac{1}{4}, \frac{1}{8}, \dots\right),$$

and

$$e^*_{2,j} \xrightarrow{\sigma(\ell_1(\Omega_0^2), W_{x_1^*, x_2^*})} x_2^* = \left(\frac{1}{2}, \frac{1}{4}, \frac{1}{8}, \dots, 0, 0, 0, \dots\right)$$

where $\left\{e^*_{1,j}, e^*_{2,j}\right\}_{j=1}^{\infty}$ is the standard basis of $\ell_1(\Omega_0^2)$. An easy computation shows that (P) holds for $\bar{x} = (1, 1, 1, \dots) \in W_{x_1^*, x_2^*}$, whereas (S) is not satisfied.

*Example* 6.8. Let $x_1^*, x_2^* \in \ell_1(\Omega_0^2)$ be defined by

$$x_1^*(i, j) = \begin{cases} \left(-\frac{1}{2}\right)^j & \text{for } i = 1 \text{ and } j \geq 1 \\ 0 & \text{for } i = 2 \text{ and } j \geq 1 \end{cases},$$

and

$$x_2^*(i, j) = \begin{cases} 0 & \text{for } i = 1 \text{ and } j \geq 1 \\ \left(-\frac{1}{3}\right)^j & \text{for } i = 2 \text{ and } j \geq 1 \end{cases}.$$

Then

$$W_{x_1^*, x_2^*} = \left\{ f \in C(\Omega^2) : f(1, 0) = \sum_{j=1}^\infty (-1)^j \frac{f\left(1, \frac{1}{j}\right)}{2^j} \text{ and } f\left(\frac{1}{2}, 0\right) = \sum_{j=1}^\infty (-1)^j \frac{f\left(\frac{1}{2}, \frac{1}{j}\right)}{3^j} \right\}.$$

By Theorem 2.1 in [5], $W_{x_1^*, x_2^*}^* = \ell_1(\Omega_0^2) = \ell_1$ and

$$e_{1,j}^* \xrightarrow{\sigma(\ell_1(\Omega_0^2), W_{x_1^*, x_2^*})} x_1^* = \left(-\frac{1}{2}, \frac{1}{4}, -\frac{1}{8}, \ldots, 0, 0, 0, \ldots\right),$$

and

$$e_{2,j}^* \xrightarrow{\sigma(\ell_1(\Omega_0^2), W_{x_1^*, x_2^*})} x_2^* = \left(0, 0, 0, \ldots, -\frac{1}{3}, \frac{1}{9}, -\frac{1}{27}, \ldots\right)$$

where $\left\{e_{1,j}^*, e_{2,j}^*\right\}_{j=1}^\infty$ is the standard basis of $\ell_1(\Omega_0^2)$. An easy computation shows that (S) holds, whereas (P) is not satisfied.

The next theorem provides a necessary and sufficient criterion for an $\ell_1$-predual space $X$ to contain a space of affine continuous functions on a Choquet simplex $S$. We also mention that this theorem can be obtained as a consequence of Theorem 3.4 in [13], but the present proof follows a different path that emphasizes the geometrical features of $W_\alpha$ spaces.

*Theorem* 6.9. An $\ell_1$-predual $X$ contains an isometric copy of a space $A(S)$ of affine continuous functions on a Choquet simplex $S$ if and only if there exist $g^* \in \left(\text{ext } B_{\ell_1}\right)'$ and $\bar{x} \in S_X$ such that $|g^*(\bar{x})| = \|g^*\| = 1$.

*Proof.* Let us suppose that $X$ contains an isometric copy of an $A(S)$ space and let $T : A(S) \to X$ be an isometrical isomorphism into. We denote by $U = T(A(S)) \subseteq X$. Since $A(S)^* = \ell_1$, there exists a subsequence $\left\{e_{n_k}^*\right\}_{k=1}^\infty$ of the standard basis of $\ell_1$ and an element $e^* \in B_{\ell_1}$ such that

$$(6.12) \qquad\qquad e_{n_k}^* \xrightarrow[k \to \infty]{\sigma(\ell_1, A(S))} e^*.$$

Let $\mathbb{1}$ denote the identically equal 1 function defined on $S$, then $e_{n_k}^*(\mathbb{1}) = 1$ for all $k \geq 1$. Thus, by (6.12), we get that $e^*(\mathbb{1}) = 1 = \|e^*\|$. Now, let us consider $u_k^* = \left(T^{-1}\right)^*\left(e_{n_k}^*\right)$ for every $k \geq 1$. It is a simple matter to show that $u_k^* \in \text{ext } B_{U^*}$ for every $k \geq 1$ and

$$(6.13) \qquad\qquad u_k^* \xrightarrow[k \to \infty]{\sigma(U^*, U)} u^*,$$

where $u^* = \left(T^{-1}\right)^* (e^*)$ and $\|u^*\| = 1$. Lemma 3.1 implies that there exist a subsequence $\left\{f^*_{n_k}\right\}^\infty_{k=1}$ of the standard basis $\{f^*_n\}^\infty_{n=1}$ of $X^* = \ell_1$ and a sequence $\{\epsilon(k)\}^\infty_{k=1}$ of signs such that $u^*_k = \epsilon(k) f^*_{n(k)}|_U$. Then, there is no loss of generality in assuming that there exists $f^* \in S_{\ell_1}$ such that

$$f^*_{n_k} \xrightarrow[k\to\infty]{\sigma(\ell_1, X)} f^*,$$

which is the desired conclusion.

Now, it remains to prove the opposite implication. If there exist $g^* \in (\text{ext } B_{\ell_1})'$ and $\bar{x} \in S_X$ such that $|g^*(\bar{x})| = \|g^*\| = 1$, then, by combining Corollary 6.5 and Remark 6.6, we have that $X$ contains an isometric copy of $W_\alpha$ and there exists $\bar{w} \in W_\alpha$ such that $\|\alpha\| = |\alpha(\bar{w})| = 1$. Remark 5.2 completes the proof. □

*Remark* 6.10. By inspecting the proof of Theorem 6.9 we obtain something more. Namely, we prove that an $\ell_1$-predual $X$ contains a space $A(S)$ of specific type whenever there exists $g^* \in (\text{ext } B_{\ell_1})'$ and $\bar{x} \in S_X$ such that $|g^*(\bar{x})| = \|g^*\| = 1$. Indeed, $X$ contains a subspace isometric to a space $W_\alpha$ where $\alpha \in S_{\ell_1}$ and it has a finite number of negative components and an infinite number of positive components or vice versa. By Remark 5.1, we can always find $\beta \in S_{\ell_1}$ such that $W_\beta$ is isometric to $W_\alpha$ and $(1, 1, 1, \dots) \in W_\beta$. Therefore, the previous theorem generalizes Theorem 1 in [38]. Indeed, the space $X = c \oplus c_0$ shows how the assumption about the existence of $g^* \in (\text{ext } B_{\ell_1})'$ and $\bar{x} \in S_X$ such that $|g^*(\bar{x})| = \|g^*\| = 1$ may be satisfied even if $\text{ext } B_X = \emptyset$.

## 7. Complemented subspaces of $\ell_1$-preduals

In this section we deal with the problem of existence of complemented copies of $W_\alpha$ spaces. We begin with the simplest situation in the setting of separable $L_1$-preduals. Let $X$ be a separable $L_1$-predual containing an isometric copy of $c$. By Proposition 3.1 in [3], there also exists a 1-complemented copy of $c$ in $X$. Moreover, we can always find a 2-complemented copy of $W_\alpha$ in $X$. Indeed, by taking into account the explicit computation of the projection constant of $W_\alpha$ spaces in $c$ (see Proposition 2.3 in [2]), there always exists an isometric copy of $W_\alpha$ such that its projection constant in $X$ does not exceed

$$1 + \left( \frac{1}{1 + \|\alpha\|} + \sum^\infty_{j=1} \frac{|\alpha(j)|}{1 + \|\alpha\| - 2|\alpha(j)|} \right)^{-1}.$$

Our next result concerns subspaces of an $\ell_1$-predual $X$ that are isometric to a specific class of $W_\alpha$ spaces (that are $M$-spaces, see Section 9) such that the norm of $\alpha$ is arbitrarily close to $r^*(X)$. Indeed, for these subspaces we do not only provide a condition guaranteeing the presence in $X$, but we also show that they are almost 1-complemented in $X$.

If $r^*(X) = 0$, we have that $X$ is isometric to a $W_\alpha$ space, where $\alpha = 0$ hence, $X$ is isometric to $c_0$. Therefore, in the sequel we assume that $r^*(X) \in (0, 1]$.

*Theorem* 7.1. Let $X$ be an $\ell_1$-predual. For every $\varepsilon \in (0, r^*(X))$ there exists a subspace $Y$ of $X$ which is isometric to $W_\alpha$ with $\alpha = \left( \frac{r^*(X)}{1+\varepsilon}, 0, 0, \dots \right)$ and $(1 + \varepsilon)$-complemented in $X$.

*Proof.* Corollary 6.3 shows that $X$ contains an isometric copy of $W_{\frac{r^*(X)}{1+\varepsilon}e_1^*}$ for every $\varepsilon \in (0, r^*(X))$. Let us consider the vectors

$$e = \left(1, \frac{r^*(X)}{1+\varepsilon}, \frac{r^*(X)}{1+\varepsilon}, \ldots, \frac{r^*(X)}{1+\varepsilon}, \ldots\right), \ e_n = (\underbrace{0, \ldots, 0}_{n-1}, 1, 0, 0, \ldots) \text{ for every } n \geq 2.$$

Clearly, $\{e\} \cup \{e_n\}_{n=2}^\infty$ forms a basis of $W_{\frac{r^*(X)}{1+\varepsilon}e_1^*}$. Let us consider the standard basis $\{e_n^*\}_{n=1}^\infty$ of $\ell_1 = W_{\frac{r^*(X)}{1+\varepsilon}e_1^*}^*$. Then, it holds

$$e_n^* \xrightarrow[n\to\infty]{\sigma\left(\ell_1, W_{\frac{r^*(X)}{1+\varepsilon}e_1^*}\right)} \frac{r^*(X)}{1+\varepsilon}e_1^*.$$

Lemma 3.1 shows that for every $e_n^*$, $n \geq 2$, there exists its extension $\widetilde{e_n^*}$ to the whole space $X$ such that $\widetilde{e_n^*} \in \text{ext } B_{X^*}$. The sequence $\left\{\widetilde{e_n^*}\right\}_{n=2}^\infty$ admits a subsequence $\left\{\widetilde{e_{n_k}^*}\right\}_{k=1}^\infty$ such that

$$\widetilde{e_{n_k}^*} \xrightarrow[k\to\infty]{\sigma(\ell_1, X)} \widetilde{e^*}.$$

Obviously, $\widetilde{e^*}$ is an extension of $\frac{r^*(X)}{1+\varepsilon}e_1^*$ and it holds

$$\frac{1}{1+\varepsilon}r^*(X) \leq \left\|\widetilde{e^*}\right\| \leq r^*(X).$$

Now, let us consider the subspace of $X$ defined as

$$Y = \overline{\text{span}\left(\{e\} \cup \{e_{n_k}\}_{k=1}^\infty\right)}.$$

Then $Y$ is isometric to $W_{\frac{r^*(X)}{1+\varepsilon}e_1^*}$. Let $P : X \to Y$ be defined by

$$P(x) = \frac{1+\varepsilon}{r^*(X)}\widetilde{e^*}(x)e + \sum_{k=1}^\infty \left(\widetilde{e_{n_k}^*}(x) - \widetilde{e^*}(x)\right)e_{n_k}.$$

for every $x \in X$. It is easily seen that the map $P$ is a projection of $X$ onto $Y$. Moreover, we have

$$\|P\| \leq \frac{1+\varepsilon}{r^*(X)}\left\|\widetilde{e^*}\right\| \leq 1+\varepsilon,$$

which completes the proof. $\qquad\square$

*Remark* 7.2. By following the proof of Theorem 7.1 with $\varepsilon = 0$, we conclude that if $X$ contains an isometric copy of $W_{r^*(X)e_1^*}$, then $X$ has a subspace $Y$ which is isometric to $W_{r^*(X)e_1^*}$ and 1-complemented in $X$.

*Remark* 7.3. Theorem 6.1 provides a sufficient condition for an $\ell_1$-predual $X$ to contain a 1-complemented isometric copy of a $W_\alpha$ space. On the other hand, Theorem 7.1 ensures that, for every $\varepsilon > 0$, there exists an isometric copy of a space $W_\beta$, with $\|\beta\| > r^*(X) - \varepsilon$, which is $(1+\varepsilon)$-complemented in $X$. Now, by taking into account the isometric embedding $T$ of $W_\alpha$ into $W_\beta$ with $\beta = \left(\frac{r^*(X)}{1+\varepsilon}, 0, 0, \ldots\right)$ considered in the proof of Proposition 5.6, we observe that $T(W_\alpha)$ is actually a hyperplane in $W_\beta$. Indeed, $T(W_\alpha) = \ker(f)$, where

$$f = \frac{r^*(X)}{1+\varepsilon}e_1^* - \sum_{i=2}^\infty \alpha(i-1)e_i^* \in \ell_1,$$

(see the duality between $W_\beta$ and $\ell_1$ described at the beginning of Section 5). Therefore, we conclude that for every $\alpha < r^*(X)$, there exists an isometric copy of $W_\alpha$ which is complemented in $X$.

## 8. $\ell_1$-PREDUALS CONTAINING ALMOST ISOMETRIC COPIES OF $c$ AND THE $w^*$-FIXED POINT PROPERTY FOR NONEXPANSIVE MAPPINGS IN $\ell_1$

A Banach space $X$ contains *almost isometric copies of* $c$ if for every $\varepsilon > 0$ there exists an isomorphic embedding $U : c \to X$ such that for every $x \in c$ we have

$$(1 - \varepsilon) \|x\| \le \|U(x)\| \le (1 + \varepsilon) \|x\| .$$

The following characterization of the presence of almost isometric copies of $c$ is a simple matter to prove.

*Proposition* 8.1. A Banach space $X$ contains almost isometric copies of $c$ if and only if, for every $\varepsilon > 0$, there exists a sequence $\{x_n\}_{n=0}^\infty$ in $X$ such that

$$(8.1) \qquad \frac{1}{1 + \varepsilon} \sup_{n \in \mathbb{N}} |t_0 + t_n| \le \left\| \sum_{n=0}^\infty t_n x_n \right\| \le \sup_{n \in \mathbb{N}} |t_0 + t_n|$$

for all $(t_n)_{n \in \mathbb{N} \cup \{0\}} \in c_0$.

A nonempty, bounded, closed and convex subset $C$ of a Banach space $X$ has the *fixed point property* (briefly, FPP) if each nonexpansive mapping $T$ (i.e., $T : C \to C$ such that $\|T(x) - T(y)\| \le \|x - y\|$ for all $x, y \in C$) has a fixed point. A dual space $X^*$ is said to have the $\sigma(X^*, X)$-*fixed point property* (briefly, $\sigma(X^*, X)$-FPP) if every nonempty, convex, $\sigma(X^*, X)$-compact set $C \subset X^*$ has the FPP.

The study of the $\sigma(X^*, X)$-FPP is especially interesting whenever the dual space has different preduals. A typical example of this situation is the space $\ell_1$ and its preduals. Indeed, $\ell_1$ has the $\sigma(\ell_1, c_0)$-FPP but it lacks the $\sigma(\ell_1, c)$-FPP. It is worth noting that a complete characterization of $\ell_1$-preduals $X$ enjoying the $\sigma(X^*, X)$-FPP has been obtained in [3] and the spaces $W_\alpha$ were deeply involved in this result. Moreover, we also consider a stronger property than the $\sigma(X^*, X)$-FPP.

*Definition* 8.2 (Definition 3.2 in [7]). A dual space $X^*$ enjoys the *stable* $\sigma(X^*, X)$-*FPP* if there exists $\gamma > 1$ such that $Y^*$ has the $\sigma(Y^*, Y)$-FPP whenever the Banach–Mazur distance $d(X, Y) < \gamma$.

Let us first recall Theorem 3.5 in [7]:

*Theorem* 8.3. Let $X$ be a predual of $\ell_1$. Then the following are equivalent.

(1) $\ell_1$ has the stable $\sigma(\ell_1, X)$-FPP.
(2) $r^*(X) < 1$.

The next theorem links the existence of almost isometric copies of $c$ in $X$ and lackness of the stable $\sigma(X^*, X)$-FPP.

*Theorem* 8.4. If a separable Banach space $X$ contains almost isometric copies of $c$, then $X^*$ lacks the stable $\sigma(X^*, X)$-FPP.

*Proof.* Fix $\varepsilon > 0$ and let $\{x_n\}_{n=0}^\infty$ in $X$ satisfy (8.1). For every $n \in \mathbb{N}$, let $x_n^* \in \overline{\left(\text{span}\left(\{x_n\}_{n=0}^\infty\right)\right)}^*$ be defined by $x_n^*(x) = t_0 + t_n$ for every $x = \sum_{n=0}^\infty t_n x_n$. From (8.1) we get $\|x_n^*\| \leq 1 + \varepsilon$. Let $\widetilde{x_n^*}$ denote a norm preserving extension of $x_n^*$ to the whole space $X$. Since $X$ is separable, there exists a subsequence $\left\{\widetilde{x_{n_k}^*}\right\}_{k=1}^\infty$ of $\left\{\widetilde{x_n^*}\right\}_{n=1}^\infty$ which is $\sigma(X^*, X)$-convergent, let us say, to $\widetilde{x_0^*}$. Clearly, $\widetilde{x_0^*}$ is an extension of $x_0^* \in \overline{\left(\text{span}\left(\{x_n\}_{n=0}^\infty\right)\right)}^*$ defined for every $x = \sum_{n=0}^\infty t_n x_n$ by $x_0^*(x) = t_0$.

Let $P : X \to \overline{\text{span}\left(\{x_0\} \cup \{x_{n_k}\}_{k=1}^\infty\right)}$ be defined by

$$P(x) = \widetilde{x_0^*}(x) x_0 + \sum_{k=1}^\infty (\widetilde{x_{n_k}^*} - \widetilde{x_0^*})(x) x_{n_k}.$$

It is easy to check that $P$ is a linear projection of $X$ onto $\overline{\text{span}\left(\{x_0\} \cup \{x_{n_k}\}_{k=1}^\infty\right)}$ with $\|P\| \leq 1 + \varepsilon$.

Let $\phi : \overline{\text{span}\left(\{x_0\} \cup \{x_{n_k}\}_{k=1}^\infty\right)} \to c$ be defined by

$$\phi\left(t_0 x_0 + \sum_{k=1}^\infty t_{n_k} x_{n_k}\right) = t_0 e_0 + \sum_{k=1}^\infty t_{n_k} e_k,$$

where

$$e_0 = (1, 1, \ldots, 1, \ldots) \text{ and } e_n = (\underbrace{0, \ldots, 0}_{n-1}, 1, 0, 0, \ldots) \text{ for every } n \in \mathbb{N}.$$

Then, $\phi$ is an isomorphism onto and for every $x \in \overline{\text{span}\left(\{x_0\} \cup \{x_{n_k}\}_{k=1}^\infty\right)}$ we have

$$\|x\| \leq \|\phi(x)\| \leq (1 + \varepsilon) \|x\|.$$

Consider now a mapping $\psi = \phi P$. Then its adjoint $\psi^* = P^* \phi^*$ is a $\sigma(\ell_1, c)$-to-$\sigma(X^*, X)$-continuous isomorphism into satisfying, for every $x \in \ell_1 = c^*$,

$$\|x\| \leq \|\psi^*(x)\| \leq (1 + \varepsilon)^2 \|x\|.$$

Let $D \subset X^*$ be defined by

$$D = \text{conv}(B_{X^*} \cup \psi^*(B_{\ell_1})).$$

It is easy to check that $D$ is convex, symmetric, $\sigma(X^*, X)$-compact, and $0$ is its interior point. Therefore, $D$ is a dual unit ball of an equivalent norm $\|\|\cdot\|\|$ on $X$. Let $Y = (X, \|\|\cdot\|\|)$, then $D = B_{Y^*}$. Clearly,

$$B_{X^*} \subset D \subset (1 + \varepsilon)^2 B_{X^*},$$

so, the Banach–Mazur distance

$$d(X, Y) \leq (1 + \varepsilon)^2.$$

Moreover, the mapping $\psi^*$ is $\sigma(\ell_1, c)$-to-$\sigma(Y^*, Y)$ continuous and since

$$B_{Y^*} \cap \psi^*(\ell_1) = \psi^*(B_{\ell_1}),$$

$\psi^*$ is an isometry from $\ell_1$ into $Y^*$.

Finally, it is well-known that a set

$$S^+ = \{x = (x(1), x(2), \dots) \in \ell_1 : \|x\| = 1 \text{ and } x(i) \geq 0 \text{ for all } i \in \mathbb{N}\}$$

is convex and $\sigma(\ell_1, c)$-compact. Moreover, a mapping $T : S^+ \to S^+$ defined by

$$T(x(1), x(2), \dots) = (0, x(1), x(2), \dots)$$

is a fixed point free isometry.

All the above shows that the set $\psi^*(S^+)$ is convex, $\sigma(Y^*, Y)$-compact and lacks the $\sigma(Y^*, Y)$-FPP, and so $X^*$ fails the stable $\sigma(X^*, X)$-FPP. □

Now we are in a position to prove the result mentioned in the introduction. It is worth noting that the assumption of Theorem 8.4 occurs to be also a necessary condition whenever $X$ is an $\ell_1$-predual.

*Theorem* 8.5. Let $X$ be a predual of $\ell_1$. Then the following are equivalent.

(1) $X$ contains almost isometric copies of $c$.
(2) $\ell_1$ lacks the stable $\sigma(\ell_1, X)$-FPP.
(3) $r^*(X) = 1$.

*Proof.* (1)$\Longrightarrow$(2) follows from Theorem 8.4. (2)$\Longleftrightarrow$(3) is a statement of Theorem 8.3. Finally, we prove the implication (3)$\Longrightarrow$(1). Let us suppose that $r^*(X) = 1$. By Corollary 6.3, the space $X$ contains an isometric copy of $W_{re_1^*}$ for every $0 \leq r < 1$. In Example 2.4 in [13], a precise value of the Banach–Mazur distance $d$ between $c$ and $W_{re_1^*}$ is given. Namely, it holds:

$$d\left(c, W_{re_1^*}\right) = \frac{3-r}{1+r}.$$

Hence, $X$ contains almost isometric copies of $c$. □

*Remark* 8.6. If we consider an $\ell_1$-predual $X$ such that $r^*(X) = 1$ and we combine Theorem 8.5 and Theorem 7.1, then we obtain that, for every $\varepsilon \in (0, 1)$, $X$ contains $(1 + \varepsilon)$-isometric copy of $c$ which is $(1 + \varepsilon)$-complemented in $X$. This fact considerably improve Theorem 2 in [38]. Indeed, there the assumption was that $B_X$ has an extreme point, which is a truly stronger requirement than the assumption $r^*(X) = 1$, as shown by the $W_\alpha$ space, where $\alpha = ((-1/2)^n)_{n=1}^\infty$. Indeed, $W_\alpha$ does not have an infinite-dimensional subspace $Y$ such that $Y^* = \ell_1$ and $\text{ext } B_Y \neq \emptyset$ (see Example 3.3 in [4]).

## 9. Classification of $L_1$-preduals

A classification of $L_1$-preduals was deeply studied in the late sixties and in the seventies (see [22], [26], [29] and the monograph [23]). These studies divided the $L_1$-preduals into several classes. Namely, the following classes were considered.

- $\mathfrak{C}$: the class of spaces $C(K)$ of all continuous real valued functions on compact Hausdorff sets $K$;
- $\mathfrak{C}_0$: the class of spaces $C_0(K)$ of all continuous real valued functions on compact Hausdorff sets $K$ which vanish at a fixed point of $K$. This class can be equivalently described by considering all the spaces of continuous real valued functions on locally compact Hausdorff spaces which vanish at infinity;

- $\mathfrak{C}_\sigma$: the class of spaces $C_\sigma(K)$ of all continuous real valued functions $f$ on compact Hausdorff sets $K$ which satisfy $f(\sigma(k)) = -f(k)$ for all $k \in K$, where $\sigma : K \to K$ is a homeomorphism of period 2 (i.e., $\sigma^2 = $ identity);
- $\mathfrak{C}_\Sigma$: the class of spaces $C_\sigma(K)$ in which the homeomorphism $\sigma$ has no fixed points;
- $\mathfrak{M}$: the class of all spaces that are isometric to sublattices of a suitable $C(K)$ space or, equivalently, the spaces $X$ which can be represented as a subspace of some $C(K)$ consisting of all functions $f$ which satisfy the following relations indexed by the elements of the set $E$:

$$(9.1) \qquad f(k_\eta^1) = \lambda_\eta f(k_\eta^2) \quad k_\eta^1, k_\eta^2 \in K, \quad \lambda_\eta \in \mathbb{R}, \ \lambda_\eta \geq 0 \quad \text{for all } \eta \in E.$$

These spaces are usually called *M*-spaces.

- $\mathfrak{G}$: the class of spaces $X$ which can be represented as a subspace of some space $C(K)$ consisting of all the functions satysfying the relations indexed by the elements of a set $E$:

$$(9.2) \qquad f(k_\eta^1) = \lambda_\eta f(k_\eta^2), \quad k_\eta^1, k_\eta^2 \in K, \quad \lambda_\eta \in \mathbb{R} \quad \text{for all } \eta \in E.$$

This class was originally introduced by Grothendieck in [15] and its elements are usually called *G*-spaces.

- $\mathfrak{A}$: the class of all the spaces $A(S)$ of affine continuous functions on compact Choquet simplexes $S$.
- $\mathfrak{A}_0$: the class of all the spaces of affine functions on compact Choquet simplexes $S$ which vanish at one fixed extreme point of $S$.

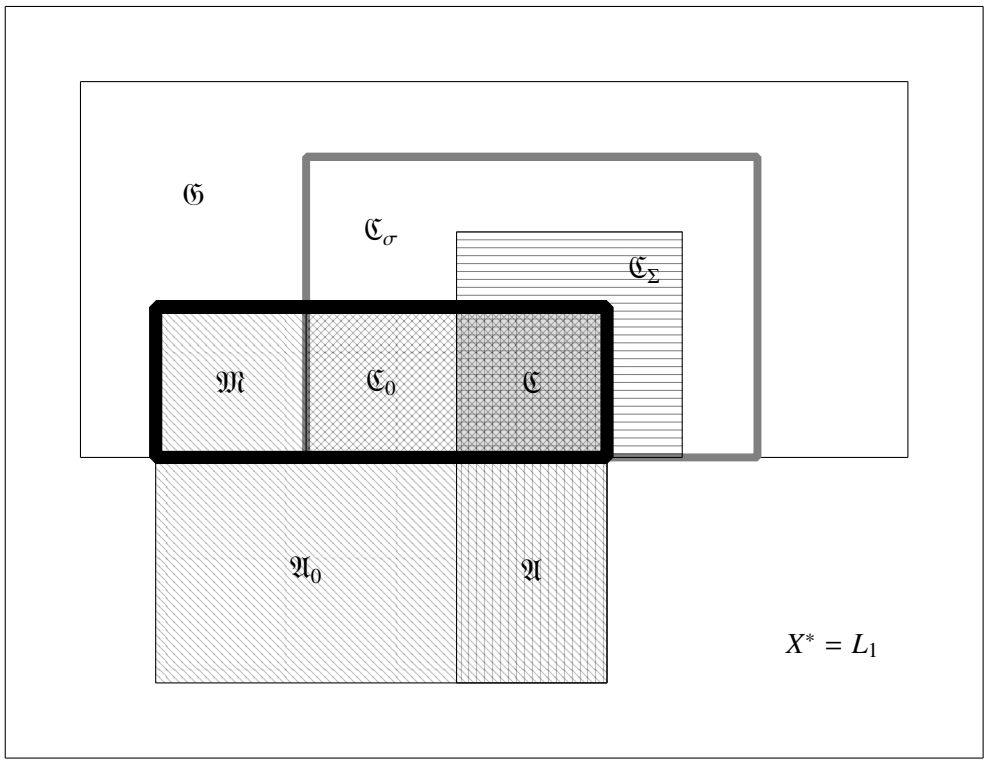

FIGURE 1. Relationships between classes of $L_1$-preduals.

In the following remark, we collect some characterizations of the classes of $L_1$-preduals listed above:

*Remark* 9.1. Let $X$ be an $L_1$-predual, then:

(i) (Corollary p. 337 in [29]) $X \in \mathfrak{C}$ if and only if ext $B_{X^*}$ is $\sigma(X^*, X)$-closed and ext $B_X \neq \emptyset$;

(ii) ([36]) $X \in \mathfrak{A}$ if and only if ext $B_X \neq \emptyset$;

(iii) ([26]) $X \in \mathfrak{A}_0$ if and only if $X^*$ is isometric to $L_1(\mu)$ in such a way that the positive cone of $L_1(\mu)$ is the image of a $\sigma(X^*, X)$-closed set.

(iv) (Theorem 1 in [29]) $X \in \mathfrak{C}_\Sigma$ if and only if ext $B_{X^*}$ is $\sigma(X^*, X)$-closed;

(v) (Theorem 13 in [11]) $X \in \mathfrak{C}_\sigma$ if and only if (ext $B_{X^*})' \subseteq$ ext $B_{X^*} \cup \{0\}$;

(vi) (Theorem 5 in [11]) $X \in \mathfrak{G}$ if and only if (ext $B_{X^*})' \subseteq \cup_{\lambda \in [0,1]} \lambda$ ext $B_{X^*}$.

The relationships between all these classes are summarized in Figure 1. The aim of this section is to show that it is possible to obtain examples of the wide variety of $L_1$-preduals even if we restrict our attention to the separable $L_1$-preduals with separable dual, i.e., the $\ell_1$-preduals. Moreover, most of these examples are spaces of type $W_\alpha$ or related to explicit models of preduals studied in [5].

We start by examining the behaviour of the spaces $W_\alpha$ with respect to the classification illustrated above. When we consider only $W_\alpha$ spaces the situation remains quite rich but a little bit simpler than the classification in the general case (see Figure 2).

*Remark* 9.2. Let us consider $\alpha \in B_{\ell_1}$. Then,

(1) $W_\alpha \in \mathfrak{C}$ if and only if $\alpha = \pm e_n^*$ for all $n \geq 1$.

(2) $W_\alpha \in \mathfrak{C}_0$ if and only if $\alpha = 0$.

(3) $W_\alpha \in \mathfrak{A}$ if and only if $\|\alpha\| = 1$ and either supp$(\alpha)$ is a finite set or the set $\mathcal{A}_{-1} = \{i \in \mathbb{N} \setminus \{0\} : \alpha(i) < 0\}$ is finite whereas $\mathcal{A}_1 = \{i \in \mathbb{N} \setminus \{0\} : \alpha(i) > 0\}$ is an infinite set (see Remark 5.2). Finally, we point out that $W_\alpha \in \mathfrak{A}$ can be seen as the space of affine continuous functions on the Choquet simplex

$$\text{(9.3)} \qquad S = \overline{\text{conv } \{\text{sgn}(\alpha(n)) e_n^*\}_{n=1}^\infty}$$

which is compact with respect to $\sigma(\ell_1, W_\alpha)$-topology.

(4) $W_\alpha \in \mathfrak{A}_0$ if and only if $\|\alpha\| \leq 1$ and either supp$(\alpha)$ is a finite set or the set $\mathcal{A}_{-1} = \{i \in \mathbb{N} \setminus \{0\} : \alpha(i) < 0\}$ is finite whereas $\mathcal{A}_1 = \{i \in \mathbb{N} \setminus \{0\} : \alpha(i) > 0\}$ is an infinite set. Indeed, let us consider the subspace of $W_\beta = A(S) \in \mathfrak{A}$ defined by

$$Y = \left\{ y = (y(1), y(2), \dots) \in W_\beta : y(n_0) = 0 \right\} = A_0(S),$$

where $S$ is the Choquet simplex defined in (9.3) and $n_0 = \min\{n \in \mathbb{N} \setminus \{0\} : \beta(n) \neq 0\}$. A trivial verification shows that $Y$ is isometric to $W_\alpha$ with $\alpha = (\beta(n))_{n=n_0+1}^\infty \in \ell_1$.

(5) $W_\alpha \in \mathfrak{M}$ if and only if $\alpha = re_n^*$ where $r \in [-1, 1]$ and $n \geq 1$. Indeed, whenever $r \geq 0$, then it follows from the definition of $\mathfrak{M}$ that $W_{re_1^*} \in \mathfrak{M}$. Moreover, if $r < 0$, to show that $W_{re_1^*} \in \mathfrak{M}$, it is sufficient to recall that Remark 5.1 proves that $W_{re_1^*}$ is isometric to $W_{|r|e_1^*}$.

It is interesting to observe that Grothendieck [15] conjectured that all the $L_1$-preduals are $G$-space, but Lindestrauss showed that this conjecture is false.

Indeed, he proved that just one of the spaces $W_\alpha$, namely those with

$$\alpha = \left(\frac{1}{2}, \frac{1}{2}, 0, 0, \dots\right),$$

is an $\ell_1$-predual (see [27], p. 81) that is not a $G$-space. Nevertheless, at least for the separable $L_1$-preduals, a weaker version of Grothendieck conjecture holds. Indeed, by Theorem 7.1, for every $\varepsilon > 0$, each separable $L_1$-predual contains a $(1 + \varepsilon)$-complemented isometric copy of $W_{\frac{r^*(X)}{1+\varepsilon}e_1^*}$, that is a $G$-space by item (5) in Remark 9.2 (more precisely, it is an $M$-space, but a $W_\alpha$ space is a $G$-space if and only if it is an $M$-space).

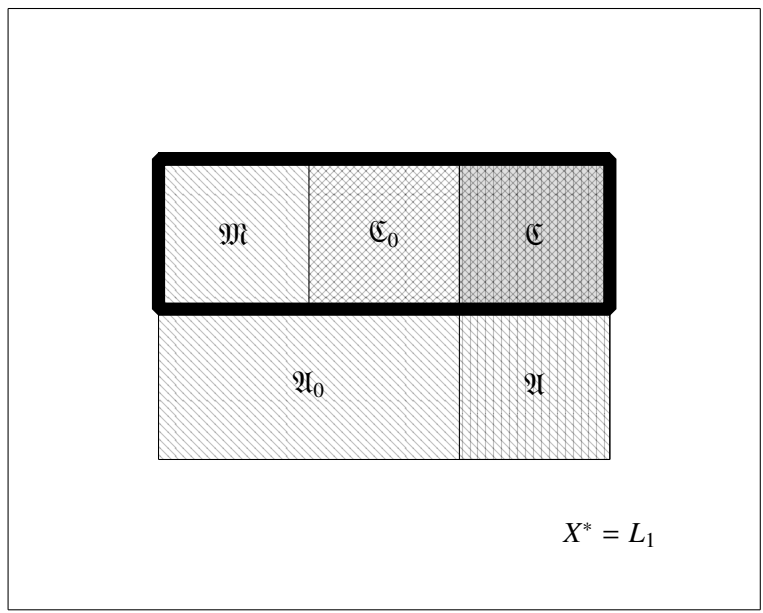

FIGURE 2. Classification of $W_\alpha$ spaces.

Remark 9.2 shows how to find examples of spaces $W_\alpha$ in the classes

$$\mathfrak{C}, \quad \mathfrak{C}_0 \setminus \mathfrak{C}, \quad \mathfrak{A} \setminus \mathfrak{C}, \quad \mathfrak{A}_0 \setminus (\mathfrak{A} \cup \mathfrak{M}), \quad \mathfrak{M} \setminus \mathfrak{C}_0.$$

Moreover, $W_\alpha \notin \mathfrak{A} \cup \mathfrak{G}$ whenever $\mathcal{A}_{-1}$ is an infinite set.

Now, we provide examples of $\ell_1$-preduals that belong to $\mathfrak{G}$ but not to $\mathfrak{C}_\sigma \cup \mathfrak{M}$.

*Example* 9.3. Let $X = W_{r_1 e_{1,1}^*, -r_2 e_{1,1}^*}$ with $0 < r_1 < r_2 < 1$ (for the notations see Example 4.2). Theorem 2.1 in [5] shows that $X^* = \ell_1$ and

$$(\text{ext } B_{X^*})' = \{r_1 e_{1,1}^*, -r_2 e_{1,1}^*\}.$$

Now, (vi) in Remark 9.1 shows that $X \in \mathfrak{G}$. Finally, we point out that $X \notin \mathfrak{M}$.

EMANUELE CASINI, Dipartimento di Scienza e Alta Tecnologia, Università dell'Insubria, via Valleggio 11, 22100 Como, Italy.

*Email address*: emanuele.casini@uninsubria.it

ENRICO MIGLIERINA, Dipartimento di Matematica per le scienze economiche, finanziarie ed attuariali, Università Cattolica del Sacro Cuore, Via Necchi 9, 20123 Milano, Italy.

*Email address*: enrico.miglierina@unicatt.it

ŁUKASZ PIASECKI, Instytut Matematyki, Uniwersytet Marii Curie-Skłodowskiej, Pl. Marii Curie-Skłodowskiej 1, 20-031 Lublin, Poland.

*Email address*: lukasz.piasecki@mail.umcs.pl